\documentclass[11pt]{article}

\usepackage{amsmath,amssymb,amsthm}
\usepackage{hyperref}
\usepackage{xcolor}
\usepackage{multirow}
\usepackage{array}
\usepackage{longtable}
\usepackage{graphicx}
\newcommand{\R}{\mathbb{R}}

\newcommand{\Om}{\Omega}
\newcommand{\B}{\mathcal{B}}
\newcommand{\Nc}{\mathcal{N}_c}
\newcommand{\vx}{\vec{x}}

\DeclareMathOperator*{\Minimize}{Minimize}

\newtheorem{teo}{Theorem}[section]
\newtheorem{lem}[teo]{Lemma}
\newtheorem{prop}[teo]{Proposition}

\theoremstyle{definition}
\newtheorem{defi}[teo]{Definition}
\newtheorem{assump}{Assumption}

\theoremstyle{remark}
\newtheorem{rem}[teo]{Remark}

\begin{document}

\title{A Telescopic Method with Guaranteed B-Stationarity\\ for Partitioned Cardinality and Complementarity Constraints\thanks{This work was
  supported by FAPESP (grants 2022/05803-3 and 2023/08706-1) and CNPq (grants 302073/2022-1
  and 305058/2024-0).}}

\author{
E. G. Birgin\thanks{Department of Computer Science, Institute of Mathematics, Statistics, and Computer Science,
  University of S\~ao Paulo, Rua do Mat\~ao, 1010, Cidade Universit\'aria, 05508-090,
  S\~ao Paulo, SP, Brazil. e-mail: egbirgin@ime.usp.br}
\and
J. M. Mart\'{\i}nez\thanks{Department of Applied Mathematics, Institute of Mathematics,
  Statistics, and Scientific Computing (IMECC), State University of Campinas, 13083-859
  Campinas SP, Brazil. e-mail: martinez@ime.unicamp.br}}

\date{October 24, 2026}

\maketitle

\begin{abstract}
We consider the minimization of a smooth function subject to box constraints and
\emph{partitioned cardinality constraints} (PCC): the variables are divided into groups and,
within each group, at most a prescribed number of variables may be non-zero.
This class contains bound-con\-strained min\-i\-mi\-za\-tion, car\-di\-nal\-ity-con\-strained min\-i\-mi\-za\-tion and
bound-con\-strained mathematical programs with complementarity constraints.
We introduce the Telescopic Method, which alternates iterations of a bound-constrained
method on the faces of the feasible set with inexpensive one-dimensional ``escaping'' tests that exchange a variable
of the current face by a currently null one within the same group.
We prove, without constraint qualifications and without requiring global solutions of
subproblems, that every limit point of the generated sequence is B-stationary, i.e., it
satisfies the first-order condition defined by the tangent cone of the feasible set.
This guarantee is obtained with exchanges of a single variable, whatever the number of
non-null entries of the limit point, because B-stationarity for this class of problems
decomposes into conditions on single variables.
On the 60 bound-constrained mathematical programs with complementarity constraints of the
test set of Kirches, Larson, Leyffer and Manns, the method found approximately B-stationary
points in all the problems.
Applications to separable regression, sparse polynomial approximation and portfolio
selection are described and tested numerically.
\end{abstract}

\section{Introduction} \label{intro}

We consider the problem
\begin{equation} \label{pcc}
\Minimize f(x) \mbox{ subject to } \|\vx_\nu\|_0 \leq p_\nu
\mbox{ and } \vec{a}_\nu \leq \vx_\nu \leq \vec{b}_\nu, \quad \nu = 1, \ldots, N,
\end{equation}
where $f: \R^n \to \R$ is continuously differentiable and $x = (\vx_1, \ldots, \vx_N)$ is a
partition of the variables into $N$ blocks (or groups)
$\vx_\nu = (x_{1,\nu}, \ldots, x_{n_\nu,\nu}) \in \R^{n_\nu}$, with
$n_1 + \cdots + n_N = n$. For each block, $\|\vx_\nu\|_0$ is the number of non-null entries
of $\vx_\nu$, $p_\nu \in \{1, \ldots, n_\nu\}$ is the maximum number of non-null variables
allowed in the block, and $\vec{a}_\nu, \vec{b}_\nu \in (\R \cup \{\pm\infty\})^{n_\nu}$ are
lower and upper bounds. We call~(\ref{pcc}) the problem with \emph{partitioned cardinality
constraints} (PCC). It contains several classical problems:
\begin{itemize}
\item \emph{Bound-constrained minimization}: $p_\nu = n_\nu$ for all $\nu$.
\item \emph{Cardinality-constrained minimization} \cite{bks}: a single constrained block
  ($p_1 < n_1$), possibly together with an unconstrained block ($p_2 = n_2$) of
  ``free'' variables.
\item \emph{Bound-constrained mathematical programs with complementarity constraints}
  \cite{kllm}: $n_\nu = 2$, $p_\nu = 1$ and $\vec{a}_\nu = 0$ for the constrained blocks,
  so that $x_{1,\nu} \geq 0$, $x_{2,\nu} \geq 0$ and $x_{1,\nu} x_{2,\nu} = 0$.
\end{itemize}

\paragraph{Motivation.}
PCC problems appear naturally in applications in which a model should keep, inside each of
several groups of terms, only a moderate number of active terms. In Section~\ref{numer} we
describe three of them: separable regression, in which each term of a sum of products of
powers of the variables may involve only a few variables; sparse polynomial approximation in
derivative-free optimization, in which the number of monomials of each degree is limited;
and portfolio selection, in which at most a prescribed number of assets of each economic
sector may be used.

\paragraph{Contributions.}
The main contributions of this paper are the following.
(i) We show that, for PCC problems, B-stationarity decomposes into conditions on single
variables (Proposition~\ref{bstat}). As a consequence, exchanges of one variable of the
current face by one variable outside it, within the same block, are enough to reach
B-stationary points, regardless of the number of non-null entries of the limit point.
(ii) We propose the \emph{Telescopic Method}, which exploits the fact that the feasible set
of~(\ref{pcc}) is a finite union of boxes (``faces''). Each exchange is evaluated by a cheap
one-dimensional procedure whose decrease is bounded below in terms of the projected
derivative of the entering variable. These tests, at most $\sum_\nu p_\nu (n_\nu - p_\nu)$
per iteration, are independent and can be performed in parallel. They are combined with a
single iteration, on each visited face, of any bound-constrained method that satisfies the
general hypotheses of Section~\ref{facemin} (``Algorithm MIA''), including Newton-type
methods such as GENCAN \cite{gencan}, which was used in our implementation.
(iii) We prove that \emph{every} limit point of the generated sequence is B-stationary,
i.e., it satisfies the first-order optimality condition defined by the tangent cone of the
feasible set, without constraint qualifications and without solving subproblems.
(iv) The same algorithm and the same theory cover cardinality-constrained problems with
arbitrary bounds, several groups with different cardinalities, and bound-constrained MPCCs.
On the 60 bound-constrained MPCCs of \cite{kllm}, the method found approximately B-stationary
points in all the problems and, using several initial points, improved the best published
objective values in 15 of them. The method was also tested on the three applications
mentioned above.

\paragraph{Related work on cardinality constraints.}
A usual approach to cardinality-constrained problems is the continuous reformulation of
Burdakov, Kanzow and Schwartz \cite{bks}, in which auxiliary variables $y \in [0,1]^n$ with
$e^T y \geq n - p$ and $x_i y_i = 0$, $i = 1, \ldots, n$, replace the constraint
$\|x\|_0 \leq p$. Several works build on this \emph{relaxed} formulation.
\v{C}ervinka, Kanzow and Schwartz \cite{cks16} studied constraint qualifications and optimality
conditions for it. Under constraint qualifications tailored to the reformulation, the limit
points of a Scholtes-type regularization method \cite{bbcs18} are S-stationary, and the
feasible limit points of an augmented Lagrangian method \cite{krsal21} are M-stationary (and
S-stationary for a suitable choice of the auxiliary variables $y$).
Krulikovski, Ribeiro and Sachine \cite{krs21} analyze constraint qualifications for the
relaxed problem (proving that Guignard's and Abadie's conditions hold when the remaining
constraints are linear) and introduce a hierarchy of weak stationarity concepts, ranging
from the weakest one up to S-stationarity. Ribeiro, Sachine and Krulikovski \cite{rsk22}
compare sequential optimality conditions (AKKT \cite{ahm11}, CAKKT \cite{ams10}, PAKKT
\cite{afss19}) for this class and propose
an approximate weak stationarity (AW-stationarity) concept that is a necessary optimality
condition without constraint qualifications. Garmanjani, Krulikovski and Ramos \cite{gkr25}
extend stationarity conditions and constraint qualifications to multiobjective problems
with cardinality constraints. These three works are concerned with optimality conditions,
not with algorithms. Medeiros, Ribeiro, Sachine and Secchin \cite{mrss25} introduce
second-order necessary conditions based on M-stationarity, stated in the original variables,
and use them to establish the global convergence of a second-order augmented Lagrangian
method. On the algorithmic side, Kreji\'c, Krulikovski and Raydan proposed
a penalty method \cite{kkr23} and an augmented Lagrangian method \cite{kkr25} for the
relaxed formulation with convex constraints, in which the subproblems are solved by the
spectral projected gradient method \cite{spg} combined with Dykstra's alternating
projections, with applications to portfolio selection and to variable selection in
constrained least squares. Augmented Lagrangian methods for problems with structured
geometric constraints, which include cardinality constraints, were analyzed in \cite{jkmw}.
Penalty decomposition methods \cite{lz13} handle a cardinality constraint by means of a
sparse copy of the variables, coupled with the original ones by a quadratic penalty, and
alternate minimizations with respect to the two sets of variables; the minimization with
respect to the sparse copy is a projection onto the cardinality-constrained set, while the
original method requires exact (global) minimizations with respect to the other variables,
which Lapucci, Levato and Sciandrone \cite{lls21} replaced by line searches along
gradient-related directions. For the minimization of a smooth function subject only to
$\|x\|_0 \leq p$, the limit points generated by both versions satisfy the conditions of Lu and
Zhang: there exists a set of $p$ indices that contains the support of the limit point and on
which the gradient vanishes \cite{lls21}. These conditions are weaker than B-stationarity when
the limit point has fewer than $p$ non-null entries (see Example~2.1 of \cite{lls21}).

The convergence guarantees of the algorithms of \cite{kkr23,kkr25} for the relaxed formulation are weaker than the one established here
for the Telescopic Method on the common class of car\-di\-nal\-ity-con\-strained problems with
bound constraints. The theory of the penalty method \cite{kkr23} requires the global
solution of each (nonconvex) penalized subproblem. For the augmented Lagrangian method,
it is shown in \cite{kkr25} that, under regularity assumptions, a limit point
$(\bar x, \bar y)$ is a first-order stationary point of the infeasibility measure and, if it
is feasible and regular, it is a first-order stationary (KKT) point of the relaxed problem. However, KKT points of the relaxed problem need not be B-stationary for the original
problem when the cardinality constraint is not active. For example, if $a < 0 < b$, then,
for every objective function $f$, the point $(\bar x, \bar y) = (0, e)$, where $e \in \R^n$ is the vector of ones, is feasible for the
relaxed problem, satisfies the linear independence constraint qualification (the only active
constraints are $x_i y_i = 0$ and $y_i \leq 1$, with gradients $(e^i, 0)$ and $(0, e^i)$),
and is a KKT point (take $-\nabla f(0)$ as the multipliers of $x_i y_i = 0$ and zero for
the other constraints). Nevertheless, if $\nabla f(0) \neq 0$, then $\bar x = 0$ is not
B-stationary for the original problem and, hence, not a local minimizer. In contrast, we
prove (Theorem~\ref{teoconv}) that every limit point generated by the Telescopic Method
is B-stationary (Proposition~\ref{bstat}), without constraint qualifications and without
solving subproblems.

\paragraph{Swap-based methods.}
Exchanging a variable of the support by a variable outside it is a natural idea for
sparsity constraints. For the minimization of a smooth function subject to
$\|x\|_0 \leq p$ (without bounds), Beck and Eldar \cite{be13} introduced coordinate-wise
minima (points at which no change of a single variable and no swap of two variables
decreases $f$) and sparse-simplex methods that perform such moves by means of exact
one-dimensional minimizations. Beck and Hallak \cite{bh16} extended this analysis to sparsity
constraints on symmetric convex sets (which include the nonnegative orthant and boxes with
the same bounds for all the variables), with methods that require the exact minimization of
$f$ on each candidate support. Swaps are also used, as a local combinatorial search, in
coordinate descent methods for $\ell_0$-regularized least squares \cite{hm20}. The Telescopic
Method combines swaps with any bound-constrained method that satisfies the hypotheses of
Section~\ref{facemin} and with one-dimensional escaping tests based on a sufficient decrease
condition. It handles several groups with different cardinalities and arbitrary bounds, and
its convergence theory requires neither exact minimizations on the supports nor exact
one-dimensional minimizations.

\paragraph{Comparison with sparse neighborhood search.}
The closest approach to ours is the framework of Lapucci, Levato, Rinaldi and
Sciandrone \cite{llrs23} for the minimization of $f(x)$ subject to $\|x\|_0 \leq s$ and
$x \in X$, with $X$ closed and convex. Their optimality condition (N-stationarity) is defined
by a neighborhood $N_\rho$ of supports at Hamming distance at most $\rho$, and their Sparse
Neighborhood Search (SNS) alternates projected gradient steps on the current support with the
exploration of the supports in $N_\rho$, each candidate being evaluated by a sequence of
projected gradient steps in the corresponding subspace. Limit points of the unsuccessful
iterations are N-stationary and, if $\rho \geq 2(s - \delta^*)$, where $\delta^*$ is the
smallest number of non-null entries of these limit points, they are basic feasible
\cite{llrs23}, a condition that, when $X$ is a box, coincides with B-stationarity. The two
approaches differ in three respects. First, the problem class: PCC involves several
cardinality constraints, one per block, and contains bound-constrained MPCCs, which are not
of the form considered in \cite{llrs23}. Second, the size of the neighborhood and the cost
of exploring it: by Proposition~\ref{bstat}, single exchanges within a block are enough for
B-stationarity at every limit point, whatever the number of its non-null entries, so the
Telescopic Method only needs at most $\sum_\nu p_\nu (n_\nu - p_\nu)$ independent
one-dimensional tests per iteration, whereas the guarantee of \cite{llrs23} requires a radius
that depends on the unknown sparsity of the limit points and a number of candidates that
grows combinatorially with it. Third, the minimization on the faces: our theory allows any
bound-constrained method that satisfies the hypotheses of Section~\ref{facemin}, of which a
single iteration is performed on each visited face, and applies to every limit point of the
sequence. Moreover, when $X$ is defined by constraints, the theory of \cite{llrs23} requires
a constraint qualification.

\paragraph{Related work on complementarity constraints.}
For mathematical programs with complementarity constraints (MPCCs) (see, e.g., \cite{lpr96}),
B-stationarity, i.e., the absence of feasible first-order descent directions, is necessary
for local optimality without any constraint qualification \cite{ss00}, whereas the
multiplier-based conditions (C-, M- and S-stationarity) are necessary only under constraint
qualifications, and C- and M-stationary points may admit feasible descent directions. Many
methods based on nonlinear programming formulations, regularizations, or relaxations of
MPCCs guarantee, under constraint qualifications, convergence to weaker stationary points,
such as C- or M-stationary points (see, e.g., \cite{atw07} and the references therein),
although standard nonlinear programming solvers are often effective in practice \cite{fl04}.
This is also the case of augmented Lagrangian methods applied to the nonlinear programming
formulation of MPCCs. Izmailov, Solodov and Uskov \cite{isu12} proved that their feasible
limit points are S-stationary under MPCC-LICQ when the multipliers remain bounded, and only
C-stationary otherwise; Andreani, Secchin and Silva \cite{ass18} replaced MPCC-LICQ by a
weaker condition (MPCC-RCPLD) in the bounded case and proved that, with unbounded
multipliers, a second-order augmented Lagrangian method converges to M-stationary points
(still under MPCC-LICQ). Andreani, Haeser, Secchin and Silva \cite{ahss19} introduced
sequential optimality conditions tailored to MPCCs (AW-, AC- and AM-stationarity), which hold
at local minimizers without constraint qualifications, together with the associated
constraint qualifications, and showed that several MPCC methods generate AC-stationary
points.

Bound-constrained mathematical programs with complementarity constraints are the PCC
problems in which every constrained block has $n_\nu = 2$, $p_\nu = 1$, $\vec{a}_\nu = 0$,
and $\vec{b}_\nu = +\infty$. For this class, Kirches, Larson, Leyffer and Manns \cite{kllm}
proposed a sequential linearization method (SLM). At each iteration, SLM solves a linear
program with complementarity constraints restricted to an $\ell_\infty$ trust region, in
order to estimate the active set. This subproblem decomposes into one- and two-dimensional
problems and is solved in linear time. The step may be improved by a bound-constrained
quadratic program on the estimated active set. Assuming only that $\nabla f$ is locally
Lipschitz continuous, it is proved in \cite{kllm} that every accumulation point of the sequence
generated by SLM is B-stationary. Since this class is a particular case of~(\ref{pcc}),
SLM and the Telescopic Method apply to the same problems in the complementarity case and
both find B-stationary points. The two methods explore the combinatorial structure in
different ways: SLM estimates the active set by means of the linear model of $f$ (the LPCC
subproblem), whereas the Telescopic Method minimizes $f$ on the current face by a
bound-constrained solver (GENCAN in our implementation) and changes the face by means of
one-dimensional escaping tests that use function and gradient evaluations. The main
difference is the scope: the decomposition of the LPCC into one- and two-dimensional
problems relies on the complementarity structure, whereas the Telescopic Method applies to
the whole PCC class, in which blocks may have any number of variables and more than one of
them may be non-null (in particular, to cardinality-constrained problems with bounds).

Other methods that compute B-stationary points of MPCCs include the $\varepsilon$-active-set
method of Fukushima and Tseng \cite{ft02} for linear complementarity constraints, under a
uniform linear independence constraint qualification; the method of Nurkanovi\'c and
Leyffer \cite{nl26}, which solves a sequence of nonlinear programs associated with branches of
the feasible set and, at each iteration, a linear program with equilibrium constraints that
either certifies B-stationarity or provides a new active set, under MPEC-MFCQ; and randomized
branch methods for linear and quadratic programs with complementarity constraints, based on
linear programming subproblems \cite{hlzz26}. Active-set methods that keep the
complementarity constraints satisfied along the iterations were also proposed by J\'udice,
Sherali, Ribeiro and Faustino~\cite{jsrf07}. In contrast with these methods, the Telescopic Method requires neither
constraint qualifications nor the solution of (linear) programs with complementarity
constraints.

In the 60 bound-constrained MPCCs of the test set of \cite{kllm} (Section~\ref{numer}),
started at the initial points of the AMPL models of the test set (with the
complementarity constraints enforced by setting one variable of each pair to zero when
necessary; see Section~\ref{numermpcc}), the Telescopic Method found approximately B-stationary points in all the problems. Since SLM and the Telescopic Method guarantee only
B-stationarity, nothing in either method makes it find better or worse B-stationary points
than the other. As an interesting fact, the objective values found by the Telescopic Method
matched (up to the two decimals reported in \cite{kllmtr}) or improved the best ones obtained
by the methods compared in \cite{kllm} in 41 of these problems, and strictly improved them in
10. Using, in addition, 20 random initial points per problem, the best values found matched or
improved the published ones in 55 of the 60 problems, and, in 15 of them, we found points with
objective values smaller than the best ones reported in \cite{kllmtr}.

In the applications, which involve groups with more than two variables and more than one
non-null variable per group, the Telescopic Method was particularly useful as a refinement of
the points given by hard thresholding or by $\ell_1$ minimization. In separable regression
(Section~\ref{numersepreg}), in the configurations in which the model without cardinality
constraints has more parameters than observations, hard thresholding failed to recover the
true support in all the instances, whereas the Telescopic Method started at the thresholding
point recovered it in 5 and 1 of the 10 instances and obtained the smallest test errors. In
sparse polynomial approximation (Section~\ref{numersparsepol}), started at the $\ell_1$
solution, the method produced polynomials that satisfy the cardinality constraints with test
errors of the same order as the ones of the $\ell_1$ solutions, which do not satisfy them.
In portfolio selection with at most $p$ assets per sector (Section~\ref{numerportfolio}),
with real data of 497 assets of the S\&P~500, the best of the solutions obtained from two
initial points was at least as good as the hard thresholding solution in all the instances,
and allowing two or three assets per sector cost only 2\% to 5\% of the risk-adjusted return
of the unconstrained portfolios. In all the applications, the runs from several initial
points showed that the difficulty is not the convergence to B-stationary points, which was
always achieved, but the existence of B-stationary points that are not global minimizers.

\medskip
The paper is organized as follows. Section~\ref{probl} states the assumptions and the
stationarity concept. Section~\ref{escape} describes the one-dimensional escaping
procedure and Section~\ref{facemin} the hypotheses on the method used to minimize on a
face. The Telescopic Method is
defined in Section~\ref{algo} and its convergence is analyzed in Section~\ref{conv}.
Section~\ref{numer} describes applications and numerical
experiments, and Section~\ref{concl} presents conclusions.

\paragraph{Notation.}
$P_{[\alpha,\beta]}(t)$ is the projection of $t \in \R$ onto the interval $[\alpha,\beta]$
(possibly unbounded). For $u, v \in (\R\cup\{\pm\infty\})^n$,
$[u, v] = \{x \in \R^n \mid u_i \leq x_i \leq v_i, \; i = 1, \ldots, n\}$.
$e^{i,\nu}$ is the canonical vector of $\R^n$ associated with the variable $x_{i,\nu}$.
$e$ is the vector of ones and $e^i$ is the $i$-th canonical vector, of appropriate dimensions.
$\#S$ is the number of elements of the finite set $S$.
$\|\cdot\|$ is the Euclidean norm.

\section{The problem and its stationarity condition} \label{probl}

In this section we state the assumptions used throughout the paper, define the projected
derivatives and the faces of the feasible set, and introduce the stationarity concept that
will be used to analyze the Telescopic Method. We show that it is a necessary optimality
condition (Proposition~\ref{locmin}) and that it coincides with B-stationarity
(Proposition~\ref{bstat}).

We denote by $\Nc = \{\nu \mid p_\nu < n_\nu\}$ the set of blocks in which the
cardinality constraint is effective. Blocks with $p_\nu = n_\nu$ are subject only to
bounds.

\begin{assump} \label{abox}
For all $\nu \in \Nc$ and $i = 1, \ldots, n_\nu$, $a_{i,\nu} \leq 0 \leq b_{i,\nu}$ and
$a_{i,\nu} < b_{i,\nu}$. For all $\nu \notin \Nc$ and $i = 1, \ldots, n_\nu$,
$a_{i,\nu} \leq b_{i,\nu}$.
\end{assump}

Assumption~\ref{abox} guarantees that setting to zero any variable of a constrained block
preserves the bound constraints. We denote by
\[
\B = [a, b] \quad \mbox{and} \quad
\Om = \{ x \in \B \mid \|\vx_\nu\|_0 \leq p_\nu, \; \nu = 1, \ldots, N\}
\]
the box and the feasible set of~(\ref{pcc}), respectively.

\begin{assump} \label{alip}
$f$ is continuously differentiable on an open set containing $\B$ and $\nabla f$ satisfies
the Lipschitz condition
$\|\nabla f(x) - \nabla f(y)\| \leq L \|x - y\|$ for all $x, y \in \B$.
\end{assump}

\begin{assump} \label{alevel}
The level set $\{x \in \Om \mid f(x) \leq f(x^0)\}$ is bounded, where $x^0 \in \Om$ is the
initial point of the algorithm.
\end{assump}

Assumption~\ref{alevel} holds trivially if $\B$ is bounded. Assumption~\ref{alip} is used
only in Lemma~\ref{lemesc}, along segments contained in $\B$.

\paragraph{Projected derivatives.}
For all $x \in \B$, $\nu = 1, \ldots, N$ and $i = 1, \ldots, n_\nu$, we define
\begin{equation} \label{gdef}
g_{i,\nu}(x) = P_{[a_{i,\nu}, b_{i,\nu}]}\left(x_{i,\nu} - \frac{\partial f}{\partial x_{i,\nu}}(x)\right) - x_{i,\nu}.
\end{equation}
By the continuity of $\nabla f$ and of the projection, $g_{i,\nu}$ is continuous on $\B$.
Moreover, $g_{i,\nu}(x) = 0$ if and only if $x_{i,\nu}$ satisfies the first-order
optimality condition for minimizing $f$ with respect to $x_{i,\nu} \in [a_{i,\nu}, b_{i,\nu}]$,
the other variables being fixed.

\begin{defi} \label{defstat}
A feasible point $x^*$ of~(\ref{pcc}) is \emph{stationary} if, for all $\nu = 1, \ldots, N$,
\begin{equation} \label{s1}
x^*_{i,\nu} \neq 0 \;\Rightarrow\; g_{i,\nu}(x^*) = 0, \quad i = 1, \ldots, n_\nu,
\end{equation}
and
\begin{equation} \label{s2}
\|\vx^*_\nu\|_0 < p_\nu \;\Rightarrow\; g_{i,\nu}(x^*) = 0 \mbox{ for all } i = 1, \ldots, n_\nu.
\end{equation}
\end{defi}

Note that, if $\nu \notin \Nc$, (\ref{s1}) and (\ref{s2}) together say that
$g_{i,\nu}(x^*) = 0$ for all $i$, i.e., the usual first-order condition with respect to
bound constraints. In the complementarity case ($n_\nu = 2$, $p_\nu = 1$,
$\vec{a}_\nu = 0$, $\vec{b}_\nu = +\infty$), if $x^*_{1,\nu} > 0$ then
(\ref{s1}) says that $\partial f(x^*)/\partial x_{1,\nu} = 0$, while if
$x^*_{1,\nu} = x^*_{2,\nu} = 0$ then (\ref{s2}) says that both partial derivatives are
non-negative. Thus, Definition~\ref{defstat} reduces to the usual B-stationarity for
complementarity constraints. Proposition~\ref{bstat} below shows that this is true in general.

\begin{prop} \label{locmin}
If $x^*$ is a local minimizer of~(\ref{pcc}), then $x^*$ is stationary.
\end{prop}

\begin{proof}
Let $\nu \in \{1, \ldots, N\}$ and $i \in \{1, \ldots, n_\nu\}$, and define
$\psi(t) = f(x^* + (t - x^*_{i,\nu}) e^{i,\nu})$, so that
$g_{i,\nu}(x^*) = P_{[a_{i,\nu},b_{i,\nu}]}(x^*_{i,\nu} - \psi'(x^*_{i,\nu})) - x^*_{i,\nu}$.

Assume that $x^*_{i,\nu} \neq 0$ and let
$T = \{ t \in [a_{i,\nu}, b_{i,\nu}] \mid |t - x^*_{i,\nu}| < |x^*_{i,\nu}| \}$.
For all $t \in T$ we have that $t \neq 0$ and, therefore, $x^* + (t - x^*_{i,\nu}) e^{i,\nu}$
has the same number of non-null entries in each block as $x^*$; since it also belongs to
$\B$, it is feasible. Since $T$ is a neighborhood of $x^*_{i,\nu}$ relative to
$[a_{i,\nu}, b_{i,\nu}]$, $x^*_{i,\nu}$ is a local minimizer of $\psi$ on
$[a_{i,\nu}, b_{i,\nu}]$ and, so, $g_{i,\nu}(x^*) = 0$. This proves~(\ref{s1}).

Assume now that $\|\vx^*_\nu\|_0 < p_\nu$. Then, for all $t \in [a_{i,\nu}, b_{i,\nu}]$,
the $\nu$-th block of $x^* + (t - x^*_{i,\nu}) e^{i,\nu}$ has at most
$\|\vx^*_\nu\|_0 + 1 \leq p_\nu$ non-null entries and the other blocks are not modified.
Thus, $x^*_{i,\nu}$ is a local minimizer of $\psi$ on $[a_{i,\nu}, b_{i,\nu}]$ and
$g_{i,\nu}(x^*) = 0$. This proves~(\ref{s2}).
\end{proof}

\paragraph{Faces.}
A \emph{face index} is an $N$-tuple $I = (I_1, \ldots, I_N)$ such that
$I_\nu \subseteq \{1, \ldots, n_\nu\}$ and $\# I_\nu = p_\nu$ for all $\nu$
(so $I_\nu = \{1, \ldots, n_\nu\}$ if $\nu \notin \Nc$). The corresponding face is the box
\begin{equation} \label{face}
F_I = \{ x \in \B \mid x_{i,\nu} = 0 \mbox{ whenever } i \notin I_\nu, \; \nu = 1, \ldots, N\}.
\end{equation}
By Assumption~\ref{abox}, $F_I \neq \emptyset$. Clearly, $\Om$ is the union of the
finitely many faces $F_I$, and every $x \in \Om$ belongs to at least one face (the faces are not disjoint: for example, $x = 0$ belongs to all of them). For
$x \in F_I$ we define $g_I(x) \in \R^n$ by
\begin{equation} \label{gI}
[g_I(x)]_{i,\nu} = g_{i,\nu}(x) \mbox{ if } i \in I_\nu, \quad
[g_I(x)]_{i,\nu} = 0 \mbox{ if } i \notin I_\nu.
\end{equation}
Since $F_I$ is a box in which the variables not in $I$ are fixed at zero, we have that
$g_I(x) = P_{F_I}(x - \nabla f(x)) - x$ for all $x \in F_I$, i.e., $g_I$ is the projected
gradient associated with the minimization of $f$ on $F_I$.

\paragraph{B-stationarity.}
Let $T_\Om(x)$ denote the (Bouligand) tangent cone of $\Om$ at $x \in \Om$. A point
$x^* \in \Om$ is said to be \emph{B-stationary} for~(\ref{pcc}) if
$\nabla f(x^*)^T d \geq 0$ for all $d \in T_\Om(x^*)$. Every local minimizer
of~(\ref{pcc}) is B-stationary.

\begin{prop} \label{bstat}
A feasible point $x^*$ of~(\ref{pcc}) is stationary in the sense of
Definition~\ref{defstat} if and only if it is B-stationary.
\end{prop}

\begin{proof}
Let $\mathcal{I}(x^*)$ be the set of face indices $I$ such that $x^* \in F_I$. Since the
faces that do not contain $x^*$ are closed and finitely many, there exists a neighborhood
$U$ of $x^*$ such that $\Om \cap U = \bigcup_{I \in \mathcal{I}(x^*)} F_I \cap U$. Hence,
$T_\Om(x^*) = \bigcup_{I \in \mathcal{I}(x^*)} T_{F_I}(x^*)$. Since each $F_I$ is a box,
$\nabla f(x^*)^T d \geq 0$ for all $d \in T_{F_I}(x^*)$ if and only if $g_I(x^*) = 0$.
Therefore, $x^*$ is B-stationary if and only if $g_{i,\nu}(x^*) = 0$ for all
$(i, \nu)$ such that $i \in I_\nu$ for some $I \in \mathcal{I}(x^*)$.

Let $S_\nu = \{ i \mid x^*_{i,\nu} \neq 0\}$. A face index $I$ belongs to
$\mathcal{I}(x^*)$ if and only if $S_\nu \subseteq I_\nu$ for all $\nu$, and the blocks
can be chosen independently. If $\# S_\nu = p_\nu$, the only possible choice is
$I_\nu = S_\nu$. If $\# S_\nu < p_\nu$, every $i \in \{1, \ldots, n_\nu\}$ belongs to some
$I_\nu \supseteq S_\nu$ with $\# I_\nu = p_\nu$. Thus, the pairs $(i, \nu)$ above are those with
$i \in S_\nu$, together with all the pairs of the blocks such that
$\|\vx^*_\nu\|_0 < p_\nu$. This is exactly what~(\ref{s1}) and~(\ref{s2}) require.
\end{proof}

Proposition~\ref{bstat} is what makes one-dimensional exchanges sufficient: since the
box constraints are separable, B-stationarity at $x^*$ reduces to conditions on single
variables, each of which can be tested by moving one variable. For a general convex set $X$,
as in \cite{llrs23}, this decomposition is not available.

\section{One-dimensional escaping procedure} \label{escape}

The Telescopic Method leaves the current face by setting to zero a variable of the face and
moving a variable outside the face away from zero. In this section we describe the
one-dimensional procedure that computes the new value of the entering variable and show that
the decrease of the objective function that it produces is bounded below in terms of the
projected derivative of that variable (Lemma~\ref{lemesc}).

Let $\alpha \leq 0 \leq \beta$ with $\alpha < \beta$ (possibly infinite), and let
$\varphi: [\alpha, \beta] \to \R$ be continuously differentiable with
\begin{equation} \label{lip1}
|\varphi'(s) - \varphi'(t)| \leq L |s - t| \mbox{ for all } s, t \in [\alpha, \beta].
\end{equation}
Then, for all $t \in [\alpha, \beta]$,
\begin{equation} \label{tay1}
\varphi(t) \leq \varphi(0) + \varphi'(0) t + \frac{L}{2} t^2.
\end{equation}
We define $g_\varphi(0) = P_{[\alpha,\beta]}(-\varphi'(0))$.

\bigskip
\noindent
{\bf Algorithm \ref{escape}.1 (Escaping)}
Let $\alpha_{esc} > 0$ and $\sigma_{\min} > 0$ be given, and assume that $g_\varphi(0) \neq 0$.

\noindent
{\bf Step 1.} Set $\sigma \leftarrow \sigma_{\min}$.

\noindent
{\bf Step 2.} Compute $t = P_{[\alpha,\beta]}(-\varphi'(0)/\sigma)$, the solution of
$\Minimize \; \varphi(0) + \varphi'(0) t + \sigma t^2/2$ subject to $t \in [\alpha, \beta]$.

\noindent
{\bf Step 3.} If
\begin{equation} \label{escdec}
\varphi(t) \leq \varphi(0) - \alpha_{esc} t^2,
\end{equation}
stop returning $t_{out} = t$. Otherwise, set $\sigma \leftarrow 2\sigma$ and go to Step~2.

\begin{lem} \label{lemesc}
Assume that~(\ref{lip1}) holds and $g_\varphi(0) \neq 0$. Then
Algorithm~\ref{escape}.1 stops after a finite number of iterations, the final value of
$\sigma$ satisfies $\sigma \leq \bar\sigma := \max\{\sigma_{\min}, L + 2\alpha_{esc}\}$, and
\begin{equation} \label{escbound}
\varphi(t_{out}) \leq \varphi(0) - \theta \, g_\varphi(0)^2,
\quad \mbox{where} \quad
\theta = \frac{\alpha_{esc}}{\max\{1, \bar\sigma\}^2}.
\end{equation}
\end{lem}

\begin{proof}
Since $g_\varphi(0) \neq 0$, we have that $\varphi'(0) \neq 0$. Assume that
$\varphi'(0) < 0$; the case $\varphi'(0) > 0$ is symmetric (with $\alpha$ in the role of
$\beta$ and negative steps). In this case $g_\varphi(0) = \min\{\beta, -\varphi'(0)\}$, so
$\beta > 0$, and for every $\sigma > 0$ the trial step is
$t = \min\{\beta, -\varphi'(0)/\sigma\} \in (0, \beta]$. Since $t \leq -\varphi'(0)/\sigma$,
we have that $\varphi'(0) t \leq -\sigma t^2$. Then, by~(\ref{tay1}),
\[
\varphi(t) \leq \varphi(0) - \left(\sigma - \frac{L}{2}\right) t^2,
\]
so~(\ref{escdec}) holds whenever $\sigma \geq L/2 + \alpha_{esc}$. Since $\sigma$ is doubled
at each failure, starting from $\sigma_{\min}$, the algorithm stops with
$\sigma \leq \max\{\sigma_{\min}, 2(L/2 + \alpha_{esc})\} = \bar\sigma$.
Moreover, if $t_{out} = \beta$ then $t_{out} \geq \min\{\beta, -\varphi'(0)\}$, and if
$t_{out} = -\varphi'(0)/\sigma$ then $t_{out} \geq \min\{\beta, -\varphi'(0)\}/\sigma$. In both
cases, $t_{out} \geq |g_\varphi(0)| / \max\{1, \bar\sigma\}$. Therefore,
by~(\ref{escdec}),
$\varphi(t_{out}) \leq \varphi(0) - \alpha_{esc} t_{out}^2 \leq \varphi(0) - \theta\, g_\varphi(0)^2$.
\end{proof}

In the Telescopic Method, Algorithm~\ref{escape}.1 is applied to
$\varphi(t) = f(\bar x + t e^{j,\nu})$ with $t \in [a_{j,\nu}, b_{j,\nu}]$, where
$\bar x \in \B$ and $\bar x_{j,\nu} = 0$. In this case,
$g_\varphi(0) = g_{j,\nu}(\bar x)$ and, by Assumption~\ref{alip},~(\ref{lip1}) holds with
the constant $L$ of Assumption~\ref{alip}.

\section{Minimization on a face: Algorithm MIA} \label{facemin}

In this section we state the hypotheses that the method used to minimize $f$ on a face
(Algorithm MIA) must satisfy.

Within a face $F_I$, problem~(\ref{pcc}) is the bound-constrained problem
\begin{equation} \label{faceprob}
\Minimize f(x) \mbox{ subject to } x \in F_I,
\end{equation}
in which the variables $x_{i,\nu}$ with $i \notin I_\nu$ are fixed at zero. Recall that,
for $x \in F_I$, $g_I(x) = P_{F_I}(x - \nabla f(x)) - x$, so that $x$ is a stationary point
of~(\ref{faceprob}) if and only if $g_I(x) = 0$.

We call \emph{Algorithm MIA} any iterative method for problems of the form~(\ref{faceprob})
that satisfies the hypotheses below. Given $w \in F_I$, an \emph{iteration} of Algorithm MIA
computes a point $\hat w \in F_I$. The iteration may depend on parameters inherited from
previous iterations (step lengths, for example); the hypotheses below must hold for all the
values of these parameters allowed by the method.

In the Telescopic Method, a single iteration of Algorithm MIA is performed on a face, the face
is abandoned, and it may be visited again later, at a point that is not the one computed by
Algorithm MIA, but one at which the objective function is not larger. For this reason, the
hypotheses concern the following kind of sequences. We say that $\{w^j\} \subseteq F_I$ is an
\emph{MIA sequence with replacements} if, for all $j$, $\hat w^{j+1} \in F_I$ is computed by
an iteration of Algorithm MIA from $w^j$ and
\begin{equation} \label{replace}
f(w^{j+1}) \leq f(\hat w^{j+1}).
\end{equation}
The choice $w^{j+1} = \hat w^{j+1}$ for all $j$ gives the sequence generated by Algorithm MIA
itself.

\begin{assump} \label{hmono}
(Descent) If $\hat w$ is computed by an iteration of Algorithm MIA from $w \in F_I$, then
$f(\hat w) \leq f(w)$.
\end{assump}

\begin{assump} \label{hconv}
(Convergence with replacements) Every limit point $w^*$ of an MIA sequence with replacements
satisfies $g_I(w^*) = 0$.
\end{assump}

\begin{assump} \label{hstep}
(Small steps) If $\{w^j\}$ is an MIA sequence with replacements and
$\lim_{j \in J} w^j = w^*$ for an infinite set $J$, then
$\lim_{j \in J} \|\hat w^{j+1} - w^j\| = 0$.
\end{assump}

Assumption~\ref{hconv} says that replacing the points computed by Algorithm MIA with points
at which $f$ is not larger does not affect its convergence. This is what allows the
Telescopic Method to perform a single iteration of Algorithm MIA on each visited face: the
point at which a face is visited again plays the role of the replacement
(Lemma~\ref{lemface}). Assumption~\ref{hconv} concerns the points from which the iterations
of Algorithm MIA start, whereas the swap search of the Telescopic Method is performed at the
points computed by these iterations. Assumption~\ref{hstep} guarantees that both converge to
the same limit (Lemma~\ref{lemstep}); it is used only there.

Assumption~\ref{hstep} is natural for bound-constrained methods. The basic iteration of a
method whose purpose is to annihilate the projected gradient is expected to produce steps
that are small when the projected gradient is small, i.e., to satisfy
\begin{equation} \label{stepbound}
\|\hat w - w\| \leq C \|g_I(w)\|
\end{equation}
for all $w \in F_I$, where $\hat w$ is computed by an iteration of Algorithm MIA from $w$ and
$C > 0$ does not depend on $w$. Under Assumption~\ref{hconv}, (\ref{stepbound}) implies
Assumption~\ref{hstep}. In fact, if $\{w^j\}$ is an MIA sequence with replacements and
$\lim_{j \in J} w^j = w^*$, then $g_I(w^*) = 0$ by Assumption~\ref{hconv}; thus, by the
continuity of $g_I$ and~(\ref{stepbound}), $\lim_{j \in J} \|\hat w^{j+1} - w^j\| \leq
C \lim_{j \in J} \|g_I(w^j)\| = 0$. Assumption~\ref{hstep} is weaker than~(\ref{stepbound})
because it concerns only the convergent subsequences of MIA sequences.

\section{The Telescopic Method} \label{algo}

In this section we define the Telescopic Method (Algorithm~\ref{algo}.1), comment on its steps
and on some variants that are useful in practice (Remark~\ref{remalgo}), and describe the
stopping criterion used in the implementation.

\noindent
{\bf Algorithm \ref{algo}.1 (Telescopic Method)}

\noindent
Let $\delta \in (0, +\infty]$, $\eta \in (0, +\infty]$ and the parameters of
Algorithm~\ref{escape}.1 be given. Let $x^0 \in \Om$ and let $I^0$ be a face index such that
$x^0 \in F_{I^0}$. Set $k \leftarrow 0$.

\medskip
\noindent
{\bf Step 1 (Swap search).}
Set $\hat x \leftarrow x^k$, $\hat I \leftarrow I^k$, and $f_{best} \leftarrow f(x^k)$.
For all $\nu \in \Nc$ and all $i \in I^k_\nu$ such that $|x^k_{i,\nu}| \leq \delta$:

\noindent
\hspace*{1em}{\bf Step 1.1.} Define $\bar x = x^k - x^k_{i,\nu} e^{i,\nu}$
(i.e., $x^k$ with the entry $x_{i,\nu}$ set to zero).

\noindent
\hspace*{1em}{\bf Step 1.2.} For all $j \in \{1, \ldots, n_\nu\} \setminus I^k_\nu$ such that
$g_{j,\nu}(\bar x) \neq 0$:

\noindent
\hspace*{2em}{\bf Step 1.2.1.} Apply Algorithm~\ref{escape}.1 to
$\varphi(t) = f(\bar x + t e^{j,\nu})$, $t \in [a_{j,\nu}, b_{j,\nu}]$, obtaining $t_{out}$,
and define $w = \bar x + t_{out} e^{j,\nu}$.

\noindent
\hspace*{2em}{\bf Step 1.2.2.} If $f(w) < f_{best}$, set $\hat x \leftarrow w$,
$f_{best} \leftarrow f(w)$, and define $\hat I$ as $I^k$ with $I^k_\nu$ replaced by
$(I^k_\nu \setminus \{i\}) \cup \{j\}$.

\noindent
\hspace*{2em}{\bf Step 1.2.3.} If $f_{best} \leq f(x^k) - \eta$, interrupt the search
and go to Step~2.

\medskip
\noindent
{\bf Step 2 (Swap).}
If $f_{best} = f(x^k)$ (no swap was found), set $z^k = x^k$ and $J^k = I^k$. Otherwise,
set $z^k = \hat x$ and $J^k = \hat I$.

\medskip
\noindent
{\bf Step 3 (Iteration of Algorithm MIA).}
Compute $x^{k+1} \in F_{J^k}$ by means of one iteration of Algorithm MIA applied
to~(\ref{faceprob}) with $I = J^k$, starting from $z^k$. Set $I^{k+1} = J^k$,
$k \leftarrow k + 1$, and go to Step~1.

\bigskip

\begin{rem} \label{remalgo}
\begin{enumerate}
\item The point $w$ of Step~1.2.1 belongs to the face $\hat I$ defined at Step~1.2.2:
  the entry $(i,\nu)$ is zero, the entry $(j,\nu)$ belongs to $[a_{j,\nu}, b_{j,\nu}]$,
  and the other entries coincide with those of $x^k \in F_{I^k}$. By Assumption~\ref{abox},
  $\bar x \in \B$. Hence all the iterates are feasible.
\item Since $\bar x_{j,\nu} = 0$, $g_\varphi(0) = g_{j,\nu}(\bar x) \neq 0$ at Step~1.2.1,
  so Algorithm~\ref{escape}.1 is well defined by Lemma~\ref{lemesc}. Moreover,
  $f(w) < f(\bar x)$. Therefore, Algorithm~\ref{algo}.1 is well defined.
\item A swap is performed only if it decreases $f$ with respect to $x^k$. Whether a swap is
  performed or not, at least one iteration of Algorithm MIA is performed on every visited
  face (Step~3). The test of Step~1.2.3 (``brute acceptance'') avoids completing the whole
  swap search when a sufficiently large decrease has already been found.
\item The results of Section~\ref{conv} remain valid, with the same proofs, if, at Step~2,
  $z^k$ is any point of a face $F_{J^k}$ such that $f(z^k) \leq f_{best}$. This allows
  the use of a \emph{look-ahead} strategy: the best candidate swaps found at Step~1 (the
  points $w$ with the smallest values of $f$) are improved by Algorithm MIA in their
  respective new faces, and $z^k$ is the best of the resulting points. Since each of these
  points has a functional value not larger than the one of the corresponding $w$, and
  $\hat x$ is among the candidates, $f(z^k) \leq f_{best}$.
\item Step~3 may perform any finite number of iterations of Algorithm MIA starting from
  $z^k$ (at least one), and $x^{k+1}$ is the last point computed. The results of
  Section~\ref{conv} remain valid with the same proofs, replacing $z^k$ with the point from
  which the last of these iterations starts.
\item All the tests of Step~1 (for different $\nu$, $i$ and $j$) are independent and can be
  performed in parallel. For each pair $(\nu, i)$, only one gradient evaluation (at $\bar x$)
  is needed to compute all the $g_{j,\nu}(\bar x)$.
\item In floating-point arithmetic, the test $g_{j,\nu}(\bar x) \neq 0$ is replaced with
  $|g_{j,\nu}(\bar x)| > \varepsilon_g$, for a small tolerance $\varepsilon_g > 0$.
\end{enumerate}
\end{rem}

\paragraph{Stopping criterion.}
Given $\varepsilon > 0$, a practical implementation stops at the first iteration $k$ in which
no swap is found at Step~1 and
\begin{equation} \label{bstatmeas}
\max \left\{ \|g_{I^k}(x^k)\|_\infty, \;
\max_{\nu \,:\, \|\vec{x}^k_\nu\|_0 < p_\nu} \; \max_{1 \leq i \leq n_\nu} |g_{i,\nu}(x^k)| \right\}
\leq \varepsilon,
\end{equation}
where $\vec{x}^k_\nu$ is the $\nu$-th block of $x^k$.
By Proposition~\ref{bstat}, the left-hand side of~(\ref{bstatmeas}) is a
B-stationarity measure: it vanishes if and only if $x^k$ is B-stationary. Note that
requiring, instead, that $|g_{j,\nu}(x^k - x^k_{i,\nu} e^{i,\nu})| \leq \varepsilon$ for
all $i \in I^k_\nu$ with $|x^k_{i,\nu}| \leq \delta$ and all $j \notin I^k_\nu$ (which
would mimic the proof of Theorem~\ref{teoconv}) is too demanding in practice: when
$x^k_{i,\nu}$ is small but not null, the corresponding swaps may fail to decrease $f$
even though $x^k$ is B-stationary.

\section{Convergence} \label{conv}

In this section we prove that every limit point of the sequence generated by the Telescopic
Method is B-stationary (Theorem~\ref{teoconv}). Throughout the section we assume that
Assumptions~\ref{abox}--\ref{alevel} hold, that Algorithm MIA satisfies
Assumptions~\ref{hmono}--\ref{hstep}, and that Algorithm~\ref{algo}.1 generates infinite
sequences $\{x^k\}$ and $\{z^k\}$.

\begin{lem} \label{lemdec}
For all $k$,
\begin{equation} \label{dec}
f(x^{k+1}) \leq f(z^k) \leq f(x^k).
\end{equation}
Consequently, $\{x^k\}$ and $\{z^k\}$ are contained in the bounded set of
Assumption~\ref{alevel}, $\{f(x^k)\}$ and $\{f(z^k)\}$ converge to the same limit, and
\begin{equation} \label{lim0}
\lim_{k \to \infty} \left[ f(x^k) - f(z^k) \right] = 0.
\end{equation}
\end{lem}

\begin{proof}
By Step~2, $f(z^k) \leq f_{best} \leq f(x^k)$ and, by Assumption~\ref{hmono},
$f(x^{k+1}) \leq f(z^k)$. Therefore, $f(z^k) \leq f(x^k) \leq f(x^0)$ for all $k$, and both
sequences lie in the level set of Assumption~\ref{alevel}, which is bounded and closed
($\Om$ is closed, as a finite union of closed sets). By continuity, $f$ is bounded below on
it, so the non-increasing sequence $\{f(x^k)\}$ converges and, by~(\ref{dec}), $\{f(z^k)\}$
converges to the same limit. This implies~(\ref{lim0}).
\end{proof}

Since there are finitely many face indices, there exist an infinite set $K$ and a face index
$I$ such that $J^k = I$ for all $k \in K$, i.e., such that the iterations of Algorithm MIA
performed at Step~3 for $k \in K$ are all performed on the same face $F_I$. The following
lemma shows that, in spite of the iterations performed on other faces between them, these
iterations behave as the iterations of Algorithm MIA applied to the minimization of $f$ on
$F_I$.

\begin{lem} \label{lemface}
Let $K$ be an infinite set and $I$ be a face index such that $J^k = I$ for all $k \in K$. Then
every limit point $x^*$ of $\{z^k\}_{k \in K}$ belongs to $F_I$ and satisfies
$g_I(x^*) = 0$.
\end{lem}

\begin{proof}
Let $k_1 < k_2 < \cdots$ be the elements of $K$, and define $w^\ell = z^{k_\ell}$ and
$\hat w^{\ell+1} = x^{k_\ell + 1}$ for all $\ell$. Then $w^\ell \in F_I$ and, by Step~3,
$\hat w^{\ell+1} \in F_I$ is computed by an iteration of Algorithm MIA from $w^\ell$. Since
$k_{\ell+1} \geq k_\ell + 1$, Lemma~\ref{lemdec} implies that
$f(w^{\ell+1}) = f(z^{k_{\ell+1}}) \leq f(x^{k_\ell + 1}) = f(\hat w^{\ell+1})$. Thus,
$\{w^\ell\}$ is an MIA sequence with replacements, and the thesis follows from
Assumption~\ref{hconv} and the fact that $F_I$ is closed.
\end{proof}

\begin{lem} \label{lemstep}
Let $K$ and $I$ be as in Lemma~\ref{lemface}, and assume that $\lim_{k \in K} z^k = x^*$.
Then $\lim_{k \in K} x^{k+1} = x^*$.
\end{lem}

\begin{proof}
With the notation of the proof of Lemma~\ref{lemface}, $\{w^\ell\}$ is an MIA sequence with
replacements that converges to $x^*$. By Assumption~\ref{hstep},
$\|x^{k_\ell+1} - z^{k_\ell}\| = \|\hat w^{\ell+1} - w^\ell\| \to 0$, and the thesis follows.
\end{proof}

\begin{teo} \label{teoconv}
The sequence $\{x^k\}$ admits at least one limit point and every limit point of $\{x^k\}$
is stationary in the sense of Definition~\ref{defstat} or, equivalently, B-stationary.
\end{teo}

\begin{proof}
By Lemma~\ref{lemdec}, $\{x^k\}$ is bounded, so it admits limit points. Let $x^*$ be one
of them, with $\lim_{k \in K_0} x^{k+1} = x^*$. Since there are finitely many face indices
and, by Lemma~\ref{lemdec}, $\{z^k\}$ is bounded, there exist an infinite set
$K \subseteq K_0$, a face index $I$ and a point $z^*$ such that $J^k = I$ for all $k \in K$
and $\lim_{k \in K} z^k = z^*$. By Lemma~\ref{lemface}, $z^* \in F_I$ and $g_I(z^*) = 0$. By
Lemma~\ref{lemstep} (this is the only place where Assumption~\ref{hstep} is used),
$\lim_{k \in K} x^{k+1} = z^*$, so $z^* = x^*$. Summing up,
\begin{equation} \label{limits}
\lim_{k \in K} z^k = \lim_{k \in K} x^{k+1} = x^* \in F_I,
\quad g_I(x^*) = 0,
\quad \mbox{and} \quad I^{k+1} = I \mbox{ for all } k \in K.
\end{equation}

\emph{Proof of~(\ref{s1}).} Let $(i,\nu)$ be such that $x^*_{i,\nu} \neq 0$. Since
$x^* \in F_I$, $i \in I_\nu$, so $g_{i,\nu}(x^*) = [g_I(x^*)]_{i,\nu} = 0$. (The same
argument shows that $g_{i,\nu}(x^*) = 0$ for all $i \in I_\nu$, $\nu = 1, \ldots, N$.)

\emph{Proof of~(\ref{s2}).} Let $\nu$ be such that $\|\vx^*_\nu\|_0 < p_\nu$. If
$\nu \notin \Nc$, then $I_\nu = \{1, \ldots, n_\nu\}$ and the thesis follows from the
observation at the end of the previous paragraph. Assume that $\nu \in \Nc$. Since
$x^*_{\ell,\nu} = 0$ for $\ell \notin I_\nu$, $\# I_\nu = p_\nu$ and
$\|\vx^*_\nu\|_0 < p_\nu$, there exists $i \in I_\nu$ such that $x^*_{i,\nu} = 0$. By the
previous paragraph, it remains to prove that $g_{\ell,\nu}(x^*) = 0$ for all
$\ell \notin I_\nu$.

Assume, by contradiction, that $\ell \notin I_\nu$ and $g_{\ell,\nu}(x^*) = c \neq 0$. For
$k \in K$, let $\bar x^k = x^{k+1} - x^{k+1}_{i,\nu} e^{i,\nu}$. Since $x^*_{i,\nu} = 0$,
by~(\ref{limits}), $\lim_{k \in K} \bar x^k = \lim_{k \in K} x^{k+1} = x^*$. By the
continuity of $f$ and $g_{\ell,\nu}$, there exists $k_0$ such that, for all $k \in K$,
$k \geq k_0$,
\[
|x^{k+1}_{i,\nu}| \leq \delta, \quad
|g_{\ell,\nu}(\bar x^k)| \geq |c|/2, \quad \mbox{and} \quad
f(\bar x^k) \leq f(x^{k+1}) + \theta c^2/8,
\]
where $\theta$ is the constant of Lemma~\ref{lemesc} with $\varphi$ as in
Section~\ref{escape}. Let $k \in K$, $k \geq k_0$. By~(\ref{limits}), $I^{k+1} = I$, so the
pair $(i, \ell)$ of block $\nu$ satisfies the conditions of Step~1 at iteration $k+1$, in
which the swap search is performed at $x^{k+1}$. Two cases are possible:
\begin{itemize}
\item The search of Step~1 is interrupted at Step~1.2.3 (possibly before examining the pair
  $(i,\ell)$). Then $f(z^{k+1}) \leq f_{best} \leq f(x^{k+1}) - \eta$.
\item The pair $(i, \ell)$ is examined. By Lemma~\ref{lemesc}, the corresponding point
  $w$ satisfies
  $f(w) \leq f(\bar x^k) - \theta (c/2)^2 \leq f(x^{k+1}) - \theta c^2/8$. Since $f_{best}$
  never increases, $f(z^{k+1}) \leq f_{best} \leq f(w)$.
\end{itemize}
In both cases, $f(x^{k+1}) - f(z^{k+1}) \geq \min\{\eta, \theta c^2/8\} > 0$ for all
$k \in K$, $k \geq k_0$, which contradicts~(\ref{lim0}). Therefore, $g_{\ell,\nu}(x^*) = 0$,
and~(\ref{s2}) is proved.
\end{proof}

The following lemma shows that, given $\varepsilon > 0$, Algorithm~\ref{algo}.1 finds, in a
finite number of iterations, a point that is approximately stationary in its face and from
which no swap decreases $f$ significantly. This is what is needed when the Telescopic
Method is used to solve subproblems approximately, as in the augmented Lagrangian approach
mentioned in Section~\ref{concl}.

\begin{lem} \label{lemfin}
For all $\varepsilon > 0$, there exists $k$ such that the swap search at iteration $k$ is not
interrupted at Step~1.2.3, $\|g_{I^k}(x^k)\|_\infty \leq \varepsilon$, and
$f(w) \geq f(x^k) - \varepsilon$ for every point $w$ computed at Step~1.2.1 of iteration $k$.
\end{lem}

\begin{proof}
If the swap search at iteration $k$ is interrupted, or if a point $w$ with
$f(w) < f(x^k) - \varepsilon$ is computed at Step~1.2.1, then, by Step~2 and
Lemma~\ref{lemdec}, $f(x^{k+1}) \leq f(z^k) \leq f(x^k) - \min\{\eta, \varepsilon\}$. Since
$\{f(x^k)\}$ is non-increasing and bounded below, this happens only for finitely many values
of $k$. Assume, by contradiction, that $\|g_{I^k}(x^k)\|_\infty > \varepsilon$ for all $k$
large enough. As in the proof of Theorem~\ref{teoconv}, there exist an infinite set $K$, a
face index $I$ and $x^*$ such that~(\ref{limits}) holds. By the continuity of $g_I$ on $F_I$,
$\lim_{k \in K} \|g_{I^{k+1}}(x^{k+1})\|_\infty = \|g_I(x^*)\|_\infty = 0$, which is a
contradiction.
\end{proof}

\begin{rem}
The threshold $\delta$ only determines which variables of the current face are candidates
to leave it, and the theory requires only $\delta > 0$. In particular, $\delta = +\infty$
(every variable of the current face is a candidate) is allowed. In our experiments
(Section~\ref{numer}), large values of $\delta$ led to much better solutions than small ones
at essentially the same cost, because swaps that remove a variable with a large absolute
value are frequently the ones that produce the largest decrease of $f$.
\end{rem}

\section{Numerical experiments} \label{numer}

In this section we report numerical experiments with the bound-constrained MPCCs of
\cite{kllm} (Section~\ref{numermpcc}) and with three applications that lead to PCC problems:
separable regression (Section~\ref{numersepreg}), sparse polynomial approximation in
derivative-free optimization (Section~\ref{numersparsepol}), and portfolio selection with
group cardinality constraints (Section~\ref{numerportfolio}). Each application is described
in the corresponding subsection, before the experiments with it. We first
describe the implementation and the parameters used in all the experiments.
All the implementations, test problems, data, and scripts used in the experiments of this
section are available at \url{https://www.ime.usp.br/~egbirgin/}.

\paragraph{Implementation.}
Algorithm~\ref{algo}.1 was implemented in Fortran 2008. Algorithm MIA is GENCAN
\cite{gencan,bmbook}, from Algencan~4.0.0~\cite{bm2020}, an active-set method for
bound-constrained minimization. Within the current face of the box, GENCAN computes Newton
directions (with line searches and extrapolations); faces are abandoned by means of spectral
projected gradient iterations \cite{spg}. Since $F_I$ is a box, GENCAN is applied
to~(\ref{faceprob}) by setting the lower and upper bounds of the variables $x_{i,\nu}$,
$i \notin I_\nu$, equal to zero (or, equivalently, by removing these variables). At Step~3,
as in Algorithm~\ref{algo}.1, a single iteration of GENCAN is performed from $z^k$ (GENCAN is
called with a limit of one iteration); one iteration of GENCAN is also performed in the
initial face before the first swap search, which amounts to replacing $x^0$ with the point so
obtained. If no improving swap is found and the stopping criterion does not hold, another
iteration is performed on the same face ($z^k = x^k$); if $f$ does not decrease in more than
three consecutive iterations on the same face, this is treated as the lack-of-progress stop of
GENCAN, which GENCAN itself would detect only after several iterations of one call. GENCAN uses second derivatives; in all the experiments the exact Hessian of $f$ was
provided, in sparse format (for the AMPL models, it is computed by the AMPL Solver Library).
Since GENCAN discards the rows and columns of the variables that are fixed at zero in the
current face, the user routine may compute the Hessian only with respect to the variables of the
face. Based on preliminary experiments, we used $\delta = +\infty$ (every variable of the current face is a
candidate to leave it) and $\eta = +\infty$ (the whole swap search is performed and the best swap
is chosen); the variant with $\delta = 0.1$ is reported below. The values
$\alpha_{esc} = 10^{-4}$ and $\sigma_{\min} = 1$ in Algorithm~\ref{escape}.1 and
$\varepsilon_g = 10^{-8}$ are usual ones. At Step~2 we use the look-ahead strategy of
Remark~\ref{remalgo} with the ten best candidates: only in the look-ahead (not at Step~3, where a
single iteration is performed), in each of the corresponding new faces,
GENCAN is run with a very large iteration limit, aiming at convergence according to the criterion
$\|g_{J}(x)\|_\infty \leq \varepsilon = 10^{-8}$, which Remark~\ref{remalgo} allows, since any $z^k$ with $f(z^k) \leq f_{best}$
may be chosen.
When the initial point has fewer than
$p_\nu$ non-null variables in a block, the initial face is completed with the null variables
with largest $|g_{i,\nu}(x^0)|$. The method stops as described in Section~\ref{algo}, with
$\varepsilon = 10^{-8}$. Besides this default version (denoted T), we report results for
variants that differ from it in a single parameter: look-ahead with a single iteration of
GENCAN in each new face (only for the MPCCs), $\delta = 0.1$, and no look-ahead. The experiments were run on a Mac Studio (model identifier Mac13,1) with an Apple M1 Max
chip (10 cores: 8 performance and 2 efficiency) and 64 GB of RAM, under macOS~13.7.8
(Darwin kernel 22.6.0, arm64). The Fortran codes (the Telescopic Method, GENCAN and the
drivers) were compiled with GNU Fortran~15.2.0 (Homebrew GCC) with option \texttt{-O3}. All runs are sequential (a
single core is used) and CPU times are in seconds.

\subsection{Bound-constrained MPCCs} \label{numermpcc}

We considered the 60 test problems of \cite{kllm}: the 40 bound-constrained QPECs (BndQPEC,
four sets of 10 randomly generated instances with 20 or 40 complementarity pairs and
indefinite or positive semidefinite Hessians) and the 20 nonlinear problems (with 20 or 40
complementarity pairs) of the BndMPCC collection. Each problem is written as a PCC problem
with one block for the variables without complementarity constraints and one block with
$n_\nu = 2$ and $p_\nu = 1$ for each complementarity pair. The initial points are taken from
the AMPL models of the test set, adapted when necessary to satisfy the bounds and the complementarity constraints: for the QPECs,
the models give no initial point and we used the AMPL default (zero), projected onto the
bounds; for the nonlinear problems, the models start with all the variables equal to one,
which violates the complementarity constraints, and we used the variables without
complementarity constraints and the first variable of each pair equal to one and the second
one equal to zero (the Rosenbrock models give no initial point and we used zero). The initial
points used by the methods compared in \cite{kllm} are not reported there; therefore, the
comparison below is with the published objective values and not with runs started from the
same points.

We compare the final objective value $f$ with the best value $f_{\mathrm{ref}}$ obtained by
the methods considered in \cite{kllm} (filterSQP, IPOPT, MINOS, SNOPT, and two variants of
SLM), as reported in Tables~5 and~9 of the technical report version of \cite{kllm}, which
contains the detailed results \cite{kllmtr}. Since these values are given with two decimals,
$f$ is rounded to two decimals before the comparison. Table~\ref{tabmpcc} reports the number
of problems in which $f_{\mathrm{ref}}$ was matched or improved, the number of problems in
which it was strictly improved, the mean and the median of the relative difference
$(f - f_{\mathrm{ref}})/\max\{1,|f_{\mathrm{ref}}|\}$, and the total CPU time over the 60
problems. Detailed results of the default version for each problem are given in
Table~\ref{tabmpccfull}. We stress that this comparison of objective function values is only
an interesting fact: SLM and the Telescopic Method both guarantee B-stationarity of the limit
points, and no property of either method favors finding better or worse B-stationary points.
Differences in the final values depend on the initial points and on the paths followed by the
methods. In all the problems the method stopped at an approximately B-stationary point.

\begin{table}[ht!]
\centering
\small
\setlength{\tabcolsep}{3.5pt}
\begin{tabular}{lccccc}
\hline
Variant & matched/improved & improved & mean diff. (\%) & median diff. (\%) & time \\
\hline
default & 41 & 10 & 0.35 & 0.00 & 0.62 \\
look-ahead with one iteration & 40 & 7 & 1.26 & 0.00 & 0.42 \\
$\delta = 0.1$ & 20 & 4 & 6.33 & 0.94 & 0.37 \\
no look-ahead & 42 & 9 & 1.63 & 0.00 & 0.28 \\
\hline
\end{tabular}
\caption{Bound-constrained MPCCs: comparison with the best values published in
\cite{kllmtr} (60 problems).}
\label{tabmpcc}
\end{table}

The default version matched or improved the best published value in 41 of the 60 problems
and strictly improved it in 10 of them; for example, in problem \texttt{rosenbrock0-40-40}
it found $f = 79.855$, whereas the best published value is $113.51$, and in
\texttt{QPEC-40-40-ind-7} it found $f = -23033$, whereas the best published value
is $-20379$. In the four Powell instances, the global minimizers (with $f = 0$) were found.
Of the 19 problems in which the best published value is smaller, 11 are QPECs with indefinite
Hessians (\texttt{ind}); in 7 of the other 8 (six QPECs with positive semidefinite Hessians
and \texttt{rosenbrock1-20-20}) the relative differences are smaller than 1\%, and in
\texttt{fletcher1-20-20} it is 2.2\%. The choice $\delta = 0.1$, that restricts the
candidates to leave the face to the variables with small absolute value, is clearly inferior.
Running GENCAN to convergence on the candidate faces of the look-ahead improves the quality of
the solutions: the mean relative difference is 0.35\% for the default version, against
1.26\% when a single iteration of GENCAN is performed in each candidate face and 1.63\%
without look-ahead, although the number of problems in which $f_{\mathrm{ref}}$ was matched
or improved is similar in the three versions; the default version found a smaller value than
these two variants in 10 and 14 problems, and a larger one in 4 and 6 problems, respectively.
The cost is moderate: the total CPU time is 0.62 seconds, against 0.42 and 0.28 seconds.
Table~\ref{tabmpccfull} shows that the method changed the face (that is, performed a swap that
decreased the objective function) in all the 40 QPECs and in 11 of the 20 nonlinear problems.
Since each outer iteration performs a single iteration of GENCAN, the number of outer
iterations (which includes the iteration performed before the first swap search) is one plus
the number of swaps plus the number of iterations in which no improving swap was found and
GENCAN was applied again on the same face. In 53 of the 60 problems, the number of outer
iterations exceeds the number of swaps by exactly one: every swap search except the last one
found an improving swap, and the last one certified approximate B-stationarity. The other 7
problems (\texttt{fletcher0-40-40} and the \texttt{mccormick} and \texttt{powell0}
instances) are nonlinear problems in which no swap was performed and between 3 and 8
iterations of GENCAN were performed on the initial face.

\paragraph{Dependence on the initial point.}
Since the method guarantees only B-stationarity, the final objective value depends on the
initial point. To quantify this dependence, each problem was also solved from 20 random
feasible initial points: each variable without complementarity constraints uniformly
distributed in the intersection of its bounds with $[v - 5, v + 5]$, where $v$ is its default
initial value, and, for each complementarity pair, one of the two variables (chosen at random)
equal to zero and the other uniformly distributed in $[0, 5]$. Table~\ref{tabmpccms}
summarizes the results. In the QPECs, the 21 runs of a problem ended, on average, at 5.42
different objective values, and a random initial point led to a value that matched or
improved the best published one in 59\% of the runs; the default initial point gave the best
of the 21 values in only 20 of the 40 QPECs. In the nonlinear problems, the default initial
point gave the best value in 18 of the 20 problems. Keeping the best of the 21 runs, the best
published value was matched or improved in 55 of the 60 problems, and strictly improved in 15
of them. Of the 1260 runs, 1259 stopped at approximately B-stationary points; the other one,
from one of the random initial points in problem \texttt{fletcher1-40-40}, whose reduced
Hessian is nearly singular \cite{kllm}, stopped because no improving swap was found although
the swap-candidates test was not satisfied (with B-stationarity measure $9.7 \times 10^{-7}$).

\begin{table}[ht!]
\centering
\small
\setlength{\tabcolsep}{4pt}
\begin{tabular}{lcccccc}
\hline
Group & problems & default & best of 21 & strictly improved & random starts & distinct values \\
\hline
QPECs & 40 & 23 & 37 & 13 & 59\% & 5.42 \\
nonlinear & 20 & 18 & 18 & 2 & 70\% & 4.05 \\
\hline
total & 60 & 41 & 55 & 15 & & \\
\hline
\end{tabular}
\caption{Bound-constrained MPCCs, default and 20 random initial points: number of problems in
which the best published value was matched or improved from the default initial point and by
the best of the 21 runs, number of problems in which the best of the 21 runs strictly improved
it, mean percentage of the random initial points from which it was matched or improved, and
mean number of different final values (rounded to four decimals) per problem.}
\label{tabmpccms}
\end{table}

{\scriptsize
\setlength{\tabcolsep}{2.5pt}
\setlength{\LTcapwidth}{\textwidth}
\begin{longtable}{rlrrrrrrrrrrr}
\caption{Bound-constrained MPCCs: detailed results of the default version of the Telescopic
Method. $n$: number of variables; $f$: final objective value; $f_{\mathrm{ref}}$: best value
published in \cite{kllmtr} (for each problem, the best of $f$, rounded to two decimals, and
$f_{\mathrm{ref}}$ is in bold; both are in bold when they coincide); B-stat.: B-stationarity
measure~(\ref{bstatmeas}) at the final point; out.: outer iterations of
Algorithm~\ref{algo}.1; GENCAN: total number of GENCAN iterations (including the look-ahead);
swaps: number of outer iterations in which the face changed (each of them decreased $f$);
\#$f$, \#$g$, \#$H$: number of evaluations of $f$, $\nabla f$, and $\nabla^2 f$; CPU: time in
seconds.}
\label{tabmpccfull}\\
\hline
\# & Problem & $n$ & $f$ & $f_{\mathrm{ref}}$ & B-stat. & out. & GENCAN & swaps & \#$f$ & \#$g$ & \#$H$ & CPU \\
\hline
\endfirsthead
\hline
\# & Problem & $n$ & $f$ & $f_{\mathrm{ref}}$ & B-stat. & out. & GENCAN & swaps & \#$f$ & \#$g$ & \#$H$ & CPU \\
\hline
\endhead
\hline
\endfoot
1 & \texttt{QPEC-20-20-ind-0} & 60 & $-6262.7792$ & $\boldsymbol{-7582.63}$ & $3.2\mathrm{e}{-14}$ & 7 & 119 & 6 & 777 & 320 & 95 & 0.01 \\
2 & \texttt{QPEC-20-20-ind-1} & 60 & $\boldsymbol{-7373.4843}$ & $\boldsymbol{-7373.48}$ & $2.8\mathrm{e}{-14}$ & 10 & 155 & 9 & 854 & 436 & 121 & 0.01 \\
3 & \texttt{QPEC-20-20-ind-2} & 60 & $\boldsymbol{-3098.2693}$ & $\boldsymbol{-3098.27}$ & $1.5\mathrm{e}{-14}$ & 12 & 156 & 11 & 1002 & 484 & 118 & 0.01 \\
4 & \texttt{QPEC-20-20-ind-3} & 60 & $\boldsymbol{-3150.7620}$ & $-3115.40$ & $7.6\mathrm{e}{-15}$ & 6 & 80 & 5 & 444 & 242 & 57 & 0.00 \\
5 & \texttt{QPEC-20-20-ind-4} & 60 & $-7038.3717$ & $\boldsymbol{-8735.77}$ & $2.1\mathrm{e}{-14}$ & 7 & 109 & 6 & 636 & 306 & 81 & 0.00 \\
6 & \texttt{QPEC-20-20-ind-5} & 60 & $-2768.9455$ & $\boldsymbol{-2829.78}$ & $3.2\mathrm{e}{-14}$ & 5 & 153 & 4 & 573 & 295 & 113 & 0.00 \\
7 & \texttt{QPEC-20-20-ind-6} & 60 & $-8172.4106$ & $\boldsymbol{-8374.58}$ & $3.6\mathrm{e}{-14}$ & 8 & 127 & 7 & 662 & 359 & 102 & 0.00 \\
8 & \texttt{QPEC-20-20-ind-7} & 60 & $-1820.6792$ & $\boldsymbol{-1923.79}$ & $1.4\mathrm{e}{-14}$ & 8 & 135 & 7 & 835 & 375 & 106 & 0.00 \\
9 & \texttt{QPEC-20-20-ind-8} & 60 & $-1967.2075$ & $\boldsymbol{-2045.32}$ & $7.1\mathrm{e}{-15}$ & 7 & 103 & 6 & 588 & 297 & 79 & 0.00 \\
10 & \texttt{QPEC-20-20-ind-9} & 60 & $-5399.2625$ & $\boldsymbol{-5622.22}$ & $1.7\mathrm{e}{-14}$ & 8 & 105 & 7 & 623 & 319 & 78 & 0.00 \\
11 & \texttt{QPEC-20-20-psd-0} & 60 & $\boldsymbol{1147.5549}$ & $\boldsymbol{1147.55}$ & $3.1\mathrm{e}{-15}$ & 8 & 41 & 7 & 485 & 239 & 33 & 0.00 \\
12 & \texttt{QPEC-20-20-psd-1} & 60 & $\boldsymbol{1043.6703}$ & $\boldsymbol{1043.67}$ & $8.9\mathrm{e}{-15}$ & 6 & 67 & 5 & 477 & 213 & 43 & 0.00 \\
13 & \texttt{QPEC-20-20-psd-2} & 60 & $\boldsymbol{1772.4854}$ & $\boldsymbol{1772.49}$ & $1.4\mathrm{e}{-14}$ & 7 & 65 & 6 & 582 & 239 & 47 & 0.00 \\
14 & \texttt{QPEC-20-20-psd-3} & 60 & $566.1266$ & $\boldsymbol{564.89}$ & $4.9\mathrm{e}{-15}$ & 13 & 124 & 12 & 1136 & 472 & 104 & 0.00 \\
15 & \texttt{QPEC-20-20-psd-4} & 60 & $\boldsymbol{896.9825}$ & $\boldsymbol{896.98}$ & $3.6\mathrm{e}{-15}$ & 6 & 40 & 5 & 424 & 182 & 26 & 0.00 \\
16 & \texttt{QPEC-20-20-psd-5} & 60 & $1489.8281$ & $\boldsymbol{1487.97}$ & $8.9\mathrm{e}{-15}$ & 7 & 47 & 6 & 446 & 217 & 38 & 0.00 \\
17 & \texttt{QPEC-20-20-psd-6} & 60 & $728.3108$ & $\boldsymbol{725.97}$ & $7.1\mathrm{e}{-15}$ & 8 & 65 & 7 & 625 & 266 & 51 & 0.00 \\
18 & \texttt{QPEC-20-20-psd-7} & 60 & $\boldsymbol{2049.9979}$ & $\boldsymbol{2050.00}$ & $7.3\mathrm{e}{-15}$ & 12 & 111 & 11 & 970 & 432 & 90 & 0.00 \\
19 & \texttt{QPEC-20-20-psd-8} & 60 & $\boldsymbol{1109.8883}$ & $\boldsymbol{1109.89}$ & $6.2\mathrm{e}{-15}$ & 8 & 36 & 7 & 399 & 229 & 31 & 0.00 \\
20 & \texttt{QPEC-20-20-psd-9} & 60 & $\boldsymbol{1075.1736}$ & $\boldsymbol{1075.17}$ & $4.4\mathrm{e}{-15}$ & 7 & 69 & 6 & 591 & 241 & 48 & 0.00 \\
21 & \texttt{QPEC-40-40-ind-0} & 120 & $-8184.0628$ & $\boldsymbol{-8254.03}$ & $1.9\mathrm{e}{-14}$ & 15 & 256 & 14 & 2268 & 999 & 204 & 0.02 \\
22 & \texttt{QPEC-40-40-ind-1} & 120 & $\boldsymbol{-18002.9442}$ & $-17970.51$ & $3.2\mathrm{e}{-14}$ & 12 & 238 & 11 & 1563 & 839 & 186 & 0.02 \\
23 & \texttt{QPEC-40-40-ind-2} & 120 & $-23315.1309$ & $\boldsymbol{-23716.62}$ & $3.7\mathrm{e}{-14}$ & 6 & 98 & 5 & 776 & 398 & 74 & 0.01 \\
24 & \texttt{QPEC-40-40-ind-3} & 120 & $\boldsymbol{-7926.9971}$ & $-7866.74$ & $2.7\mathrm{e}{-14}$ & 15 & 264 & 14 & 2207 & 1020 & 218 & 0.02 \\
25 & \texttt{QPEC-40-40-ind-4} & 120 & $\boldsymbol{-17463.0028}$ & $-17301.43$ & $2.7\mathrm{e}{-14}$ & 17 & 357 & 16 & 2707 & 1241 & 285 & 0.02 \\
26 & \texttt{QPEC-40-40-ind-5} & 120 & $\boldsymbol{-10721.1017}$ & $-10672.87$ & $1.4\mathrm{e}{-14}$ & 17 & 374 & 16 & 2520 & 1245 & 307 & 0.02 \\
27 & \texttt{QPEC-40-40-ind-6} & 120 & $-4828.5521$ & $\boldsymbol{-5073.21}$ & $3.6\mathrm{e}{-14}$ & 18 & 314 & 17 & 2874 & 1227 & 262 & 0.03 \\
28 & \texttt{QPEC-40-40-ind-7} & 120 & $\boldsymbol{-23033.4767}$ & $-20378.55$ & $7.1\mathrm{e}{-14}$ & 21 & 503 & 20 & 3646 & 1575 & 380 & 0.04 \\
29 & \texttt{QPEC-40-40-ind-8} & 120 & $-14752.6274$ & $\boldsymbol{-15355.83}$ & $5.0\mathrm{e}{-14}$ & 21 & 356 & 20 & 3207 & 1426 & 292 & 0.03 \\
30 & \texttt{QPEC-40-40-ind-9} & 120 & $\boldsymbol{-11058.4014}$ & $-10848.63$ & $2.1\mathrm{e}{-14}$ & 19 & 336 & 18 & 2922 & 1305 & 273 & 0.03 \\
31 & \texttt{QPEC-40-40-psd-0} & 120 & $4002.4672$ & $\boldsymbol{4002.46}$ & $1.9\mathrm{e}{-14}$ & 15 & 164 & 14 & 2242 & 874 & 139 & 0.02 \\
32 & \texttt{QPEC-40-40-psd-1} & 120 & $\boldsymbol{4130.5093}$ & $\boldsymbol{4130.51}$ & $1.9\mathrm{e}{-14}$ & 10 & 98 & 9 & 1271 & 555 & 76 & 0.01 \\
33 & \texttt{QPEC-40-40-psd-2} & 120 & $\boldsymbol{3105.8480}$ & $\boldsymbol{3105.85}$ & $3.6\mathrm{e}{-14}$ & 8 & 56 & 7 & 1062 & 420 & 44 & 0.01 \\
34 & \texttt{QPEC-40-40-psd-3} & 120 & $\boldsymbol{4941.1088}$ & $\boldsymbol{4941.11}$ & $1.8\mathrm{e}{-14}$ & 14 & 152 & 13 & 1906 & 822 & 127 & 0.02 \\
35 & \texttt{QPEC-40-40-psd-4} & 120 & $\boldsymbol{2446.0775}$ & $\boldsymbol{2446.08}$ & $8.9\mathrm{e}{-15}$ & 17 & 180 & 16 & 2650 & 992 & 157 & 0.02 \\
36 & \texttt{QPEC-40-40-psd-5} & 120 & $\boldsymbol{2364.7669}$ & $\boldsymbol{2364.77}$ & $2.4\mathrm{e}{-14}$ & 20 & 243 & 19 & 3027 & 1208 & 194 & 0.03 \\
37 & \texttt{QPEC-40-40-psd-6} & 120 & $\boldsymbol{4035.2263}$ & $\boldsymbol{4035.23}$ & $1.5\mathrm{e}{-14}$ & 12 & 92 & 11 & 1579 & 650 & 80 & 0.01 \\
38 & \texttt{QPEC-40-40-psd-7} & 120 & $\boldsymbol{3152.8609}$ & $3152.99$ & $1.7\mathrm{e}{-14}$ & 17 & 214 & 16 & 2933 & 1027 & 174 & 0.03 \\
39 & \texttt{QPEC-40-40-psd-8} & 120 & $2220.3469$ & $\boldsymbol{2220.32}$ & $1.2\mathrm{e}{-14}$ & 15 & 152 & 14 & 2231 & 865 & 130 & 0.02 \\
40 & \texttt{QPEC-40-40-psd-9} & 120 & $3657.8812$ & $\boldsymbol{3657.14}$ & $1.1\mathrm{e}{-14}$ & 16 & 191 & 15 & 2411 & 954 & 156 & 0.02 \\
41 & \texttt{fletcher0-20-20} & 60 & $\boldsymbol{2246.8641}$ & $\boldsymbol{2246.86}$ & $1.7\mathrm{e}{-10}$ & 2 & 459 & 1 & 901 & 502 & 459 & 0.01 \\
42 & \texttt{fletcher1-20-20} & 60 & $4017.7166$ & $\boldsymbol{3930.27}$ & $9.0\mathrm{e}{-9}$ & 4 & 631 & 3 & 1635 & 806 & 491 & 0.01 \\
43 & \texttt{himmelblau0-20-20} & 60 & $\boldsymbol{826.2919}$ & $\boldsymbol{826.29}$ & $1.2\mathrm{e}{-14}$ & 11 & 366 & 10 & 1525 & 703 & 366 & 0.01 \\
44 & \texttt{himmelblau1-20-20} & 60 & $\boldsymbol{230.3125}$ & $\boldsymbol{230.31}$ & $1.1\mathrm{e}{-13}$ & 3 & 61 & 2 & 377 & 136 & 61 & 0.00 \\
45 & \texttt{mccormick0-20-20} & 60 & $\boldsymbol{58.9346}$ & $\boldsymbol{58.93}$ & $5.9\mathrm{e}{-13}$ & 8 & 8 & 0 & 190 & 179 & 8 & 0.00 \\
46 & \texttt{mccormick1-20-20} & 60 & $\boldsymbol{58.9346}$ & $\boldsymbol{58.93}$ & $1.1\mathrm{e}{-16}$ & 5 & 5 & 0 & 119 & 111 & 4 & 0.00 \\
47 & \texttt{powell0-20-20} & 60 & $\boldsymbol{0}$ & $\boldsymbol{0}$ & $0$ & 3 & 3 & 0 & 69 & 67 & 3 & 0.00 \\
48 & \texttt{powell1-20-20} & 60 & $\boldsymbol{0}$ & $\boldsymbol{0}$ & $0$ & 1 & 1 & 0 & 26 & 23 & 1 & 0.00 \\
49 & \texttt{rosenbrock0-20-20} & 60 & $\boldsymbol{40.0567}$ & $\boldsymbol{40.06}$ & $7.8\mathrm{e}{-10}$ & 21 & 106 & 20 & 3935 & 567 & 104 & 0.00 \\
50 & \texttt{rosenbrock1-20-20} & 60 & $58.3911$ & $\boldsymbol{58.34}$ & $7.9\mathrm{e}{-10}$ & 2 & 7 & 1 & 369 & 50 & 5 & 0.00 \\
51 & \texttt{fletcher0-40-40} & 120 & $\boldsymbol{4246.8641}$ & $\boldsymbol{4246.86}$ & $2.1\mathrm{e}{-12}$ & 6 & 6 & 0 & 362 & 252 & 6 & 0.00 \\
52 & \texttt{fletcher1-40-40} & 120 & $\boldsymbol{7784.5363}$ & $\boldsymbol{7784.54}$ & $7.4\mathrm{e}{-9}$ & 20 & 1457 & 19 & 3388 & 2443 & 1457 & 0.07 \\
53 & \texttt{himmelblau0-40-40} & 120 & $\boldsymbol{1652.5837}$ & $\boldsymbol{1652.58}$ & $1.2\mathrm{e}{-14}$ & 21 & 766 & 20 & 4825 & 1833 & 766 & 0.03 \\
54 & \texttt{himmelblau1-40-40} & 120 & $\boldsymbol{460.6250}$ & $460.63$ & $1.1\mathrm{e}{-13}$ & 3 & 61 & 2 & 677 & 196 & 61 & 0.00 \\
55 & \texttt{mccormick0-40-40} & 120 & $\boldsymbol{118.9346}$ & $\boldsymbol{118.93}$ & $1.1\mathrm{e}{-16}$ & 8 & 8 & 0 & 347 & 338 & 7 & 0.00 \\
56 & \texttt{mccormick1-40-40} & 120 & $\boldsymbol{118.9346}$ & $\boldsymbol{118.93}$ & $1.1\mathrm{e}{-16}$ & 5 & 5 & 0 & 219 & 211 & 4 & 0.00 \\
57 & \texttt{powell0-40-40} & 120 & $\boldsymbol{0}$ & $\boldsymbol{0}$ & $0$ & 3 & 3 & 0 & 129 & 127 & 3 & 0.00 \\
58 & \texttt{powell1-40-40} & 120 & $\boldsymbol{0}$ & $\boldsymbol{0}$ & $0$ & 1 & 1 & 0 & 46 & 43 & 1 & 0.00 \\
59 & \texttt{rosenbrock0-40-40} & 120 & $\boldsymbol{79.8547}$ & $113.51$ & $7.8\mathrm{e}{-10}$ & 41 & 206 & 40 & 15054 & 1927 & 204 & 0.02 \\
60 & \texttt{rosenbrock1-40-40} & 120 & $\boldsymbol{117.9910}$ & $\boldsymbol{117.99}$ & $7.9\mathrm{e}{-10}$ & 2 & 7 & 1 & 729 & 90 & 5 & 0.00 \\
\hline
\multicolumn{6}{l}{Total} & 635 & 10712 & 544 & 94053 & 36409 & 9235 & 0.62 \\
\end{longtable}
}

\subsection{Separable regression} \label{numersepreg}

In many practical problems we wish to express a variable $y$ as a function of $d$ variables
$\xi_1, \ldots, \xi_d$ in the form
\begin{equation} \label{apr}
y \approx \sum_{\nu=1}^{M} c_\nu \, \xi_1^{\alpha_{1,\nu}} \cdots \xi_d^{\alpha_{d,\nu}},
\end{equation}
given a table of observations $(\xi_{\ell,1}, \ldots, \xi_{\ell,d}, y_\ell)$,
$\ell = 1, \ldots, m$, with $\xi_{\ell,j} > 0$ for all $\ell, j$. The unknowns are the linear
coefficients $c_\nu$ and the exponents $\alpha_{j,\nu}$, and the fitting problem is
\[
\Minimize \sum_{\ell=1}^m \left[ \sum_{\nu=1}^M c_\nu \prod_{j=1}^d \xi_{\ell,j}^{\alpha_{j,\nu}} - y_\ell \right]^2 .
\]
If $d$ is large, each term of~(\ref{apr}) should involve a moderate number of variables.
Therefore, for each $\nu = 1, \ldots, M$, we impose that at most $p_\nu \leq d$ of the
exponents $\alpha_{1,\nu}, \ldots, \alpha_{d,\nu}$ are non-null. This is a PCC problem with
$N = M + 1$ blocks: the blocks $\vec{\alpha}_\nu \in \R^d$, $\nu = 1, \ldots, M$, with
cardinality $p_\nu$, and the block $c \in \R^M$ without cardinality constraint. Since
$\xi^{\alpha} = e^{\alpha \ln \xi}$, the objective function is $C^\infty$. Bounds
$-\alpha_{\max} \leq \alpha_{j,\nu} \leq \alpha_{\max}$ and $-c_{\max} \leq c_\nu \leq c_{\max}$
make $\B$ bounded, so that Assumptions~\ref{alip} and~\ref{alevel} hold.

We generated instances of this problem as follows. The data
$\xi_{\ell,j}$ are uniformly distributed in $[0.5, 2]$. A ``true'' model with $M$ terms is
generated choosing, for each term, $p$ variables at random, with exponents of random sign
and absolute value uniformly distributed in $[0.5, 2]$, and coefficients $c_\nu$ of random sign
and absolute value uniformly distributed in $[1, 3]$. The observations $y_\ell$ are the
values of the true model perturbed by Gaussian noise with standard deviation equal to 1\% or
5\% of the standard deviation of the noiseless values. Four configurations were
considered, with 10 instances each: (A) $d = 10$, $M = 3$, $p = 2$, $m = 100$, noise 1\%;
(B) $d = 20$, $M = 4$, $p = 3$, $m = 200$, noise 1\%; (C) $d = 30$, $M = 4$, $p = 3$,
$m = 80$, noise 5\%; and (D) $d = 50$, $M = 5$, $p = 3$, $m = 100$, noise 5\%. In C and D the number of parameters of the model without
cardinality constraints, $M(d + 1)$, exceeds the number of observations $m$. The bounds are
$|\alpha_{j,\nu}| \leq 3$ and $|c_\nu| \leq 10$. Each method was run from five random initial
points (with $p$ random non-null exponents of absolute value at most 1 in each term), and the
solution with the smallest objective value was kept.

Besides the three versions of the Telescopic Method, we considered a hard thresholding
strategy (HT) \cite{bd2009}: the problem is solved without cardinality constraints (i.e.,
with $p_\nu = d$) by GENCAN, the $p$ exponents of largest absolute value of each term are kept,
and the objective function is minimized (by GENCAN) on the corresponding face. HT+T denotes
the Telescopic Method started at the point given by HT. Table~\ref{tabsepreg} shows the
number of instances in which the support of the true model was exactly recovered (up to a
permutation of the terms), the mean fraction of the true non-null exponents that were
identified, the median of the relative error on 2000 test points (root mean square error of
the prediction of the noiseless values divided by their standard deviation), and the mean
CPU time per instance (five starting points).

\begin{table}[ht!]
\centering
\small
\begin{tabular}{llcccc}
\hline
Config. & Method & supports recovered & fraction recovered & test error & time \\
\hline
\multirow{5}{*}{A} & T & 10 & 1.00 & 0.0029 & 0.05 \\
 & T ($\delta=0.1$) & 5 & 0.83 & 0.0765 & 0.02 \\
 & T (no l.-a.) & 9 & 0.97 & 0.0029 & 0.08 \\
 & HT & 10 & 1.00 & 0.0029 & 0.02 \\
 & HT+T & 10 & 1.00 & 0.0029 & 0.03 \\
\hline
\multirow{5}{*}{B} & T & 10 & 1.00 & 0.0031 & 0.84 \\
 & T ($\delta=0.1$) & 1 & 0.70 & 0.2846 & 0.26 \\
 & T (no l.-a.) & 8 & 0.97 & 0.0032 & 1.65 \\
 & HT & 10 & 1.00 & 0.0031 & 0.24 \\
 & HT+T & 10 & 1.00 & 0.0031 & 0.51 \\
\hline
\multirow{5}{*}{C} & T & 1 & 0.75 & 0.2712 & 0.53 \\
 & T ($\delta=0.1$) & 0 & 0.47 & 0.6024 & 0.12 \\
 & T (no l.-a.) & 1 & 0.71 & 0.3229 & 1.18 \\
 & HT & 0 & 0.59 & 0.4239 & 0.93 \\
 & HT+T & 5 & 0.92 & 0.0539 & 1.27 \\
\hline
\multirow{5}{*}{D} & T & 2 & 0.77 & 0.2381 & 2.57 \\
 & T ($\delta=0.1$) & 0 & 0.44 & 0.5292 & 0.30 \\
 & T (no l.-a.) & 0 & 0.69 & 0.2665 & 6.03 \\
 & HT & 0 & 0.43 & 0.5559 & 0.98 \\
 & HT+T & 1 & 0.83 & 0.1855 & 2.97 \\
\hline
\end{tabular}
\caption{Separable regression: 10 instances per configuration.}
\label{tabsepreg}
\end{table}

In the easy configurations A and B, HT recovers the true model in all the instances, and so
does the Telescopic Method. In the configurations C and
D, in which the model without cardinality constraints has more parameters than observations
and the noise is larger, HT fails in all the instances, while the Telescopic Method started
from the HT point (HT+T) recovers the true support in 5 and 1 of the 10 instances and has
the smallest test errors. Thus, the swaps of the Telescopic Method are able to correct
the support selected by thresholding. The variant with $\delta = 0.1$ is clearly inferior to
the default version, and so is, to a lesser extent, the variant without look-ahead, which was
also slower in all the configurations. In all the runs of the Telescopic Method (all
the versions and HT+T, from the five initial points), the method stopped at an approximately
B-stationary point, with the B-stationarity measure~(\ref{bstatmeas}) at most $10^{-8}$.

\paragraph{Why thresholding fails.}
To understand these results, we also computed, for each instance, the model without
cardinality constraints (best of the same five starts) and the least squares fit on the true
support (computed by GENCAN started at the true model). The model without cardinality
constraints uses all the $M d$ exponents. In configurations A and B it is still a reasonable
model (median test errors 0.008 and 0.012), although worse than the fit on the true support
(0.0029 and 0.0031). In configurations C and D, it interpolates the noisy observations (its
residual is about $10^{-20}$ times the one of the fit on the true support) and is useless as a
predictor (median test errors 0.97 and 0.90, against 0.026 for the fit on the true support).
Thresholding the exponents of an interpolating model selects essentially arbitrary supports,
which explains the failure of HT in C and D. In these configurations, the cardinality
constraints are what makes the problem well posed.

\paragraph{Local versus global solutions.}
In every instance in which T, HT, or HT+T did not recover the true support, the sum of squares at
its solution was larger than the one of the fit on the true support. Thus, the failures are
not due to the model (the true support is the best one among the computed solutions) but to
the fact that the methods stop at B-stationary points that are not global minimizers. This is
also reflected by the individual starts: the fraction of the starting points from which T
recovered the true support was 66\%, 34\%, 2\%, and 4\% in configurations A, B, C, and D,
respectively, and it was 100\%, 80\%, 14\%, and 2\% for HT+T. Hence, a good initial point (as
the HT point) and several starts are important in the harder configurations.

\subsection{Sparse polynomial approximation} \label{numersparsepol}

In derivative-free optimization, when the objective function $F: \R^d \to \R$ is expensive to
evaluate, it is common to use polynomial surrogates. Given $m$ previous evaluations
$F(z^1), \ldots, F(z^m)$, we construct a polynomial of degree $D \in \{2, 3, 4, 5\}$,
\[
Q(z) = \sum_{|\beta| \leq D} c_\beta \, z^\beta,
\]
and choose a candidate point by minimizing $Q$ on a trust region. We suggest limiting the
number of terms of each degree: for $s = 2, \ldots, D$, at most $p_s$ coefficients
$c_\beta$ with $|\beta| = s$ may be non-null. The coefficients are computed by solving
\[
\Minimize \sum_{\ell=1}^m \left[ \sum_{|\beta| \leq D} c_\beta (z^\ell)^\beta - F(z^\ell) \right]^2
\]
subject to these constraints. This is a PCC problem with one block per degree $s \geq 2$
(with cardinality $p_s$) and one block (without cardinality constraint) containing the
coefficients of degree $0$ and $1$. The objective function is a convex quadratic, so
Assumption~\ref{alip} holds with $L = 2\|A\|^2$, where $A$ is the $m \times \#\{\beta\}$
matrix with entries $(z^\ell)^\beta$. Since $A$ may have a nontrivial null space (in
particular, when $m$ is smaller than the number of monomials), we impose bounds
$-c_{\max} \leq c_\beta \leq c_{\max}$, so that $\B$ is bounded and Assumption~\ref{alevel}
holds.

In the experiments, we considered $d = 6$ and degree $D = 4$ (210
monomials). The sample points are $z^\ell = \bar z + \Delta s^\ell$, with $\Delta = 0.5$, $\bar z$
uniformly distributed in $[-1,1]^6$, and $s^\ell$ uniformly distributed in $[-1,1]^6$,
$\ell = 1, \ldots, m$, with $m \in \{40, 60, 90\}$, so that the number of sample points is
smaller than the number of monomials. The polynomial is written in the variables $s$, the
function values are divided by their largest absolute value, and $c_{\max} = 10^4$. Two functions were considered:
\begin{itemize}
\item ROS: the chained Rosenbrock function
  $\sum_{i=1}^{5} [100 (x_{i+1} - x_i^2)^2 + (1 - x_i)^2]$, which is a polynomial of degree
  4 with 33 non-null coefficients in the variables $s$ (11, 10, and 5 of degrees 2, 3,
  and 4, respectively); we used these numbers as $p_2$, $p_3$, and $p_4$;
\item EXP: $\sum_{i=1}^{5} \exp(0.5\, x_i x_{i+1}) + \sum_{i=1}^{6} \sin(x_i)$, which is not
  a polynomial; we used as $p_k$ the number of monomials of degree $k$ that involve a single
  variable or two consecutive variables ($p_2 = 11$, $p_3 = 16$, $p_4 = 21$).
\end{itemize}
For each function and each $m$ we generated 10 instances. We compared the Telescopic Method
with two methods that do not enforce the cardinality constraints, namely, the minimum-norm least
squares solution (MN) and the solution of the $\ell_1$ minimization problem
$\min \sum_{|\beta| \geq 2} |c_\beta|$ subject to the interpolation conditions, as proposed by
Bandeira, Scheinberg, and Vicente \cite{bsv2012} for the computation of sparse models in
derivative-free optimization (solved by the HiGHS linear programming solver \cite{hh2018}),
and with hard thresholding of MN followed by a least squares fit on the selected support (HT).
The Telescopic Method was started at $c = 0$ (T and its variants), at the HT point (HT+T), and
at the $\ell_1$ solution (L1+T). Table~\ref{tabsparsepol} shows the median, over the 10
instances, of the relative error $\|Q - F\| / \|F - \bar F\|$ on 2000 test points (where
$\bar F$ is the mean value of $F$ on the test points) and, for ROS, the number of instances in
which the true polynomial was recovered (maximum coefficient error below $10^{-5}$).

\begin{table}[ht!]
\centering
\scriptsize
\setlength{\tabcolsep}{3pt}
\begin{tabular}{llcccccccc}
\hline
$F$ & $m$ & MN & L1 & HT & T & T ($\delta = 0.1$) & T (no l.-a.) & HT+T & L1+T \\
\hline
ROS & 40 & $2.4\mathrm{e}{-1}$ (0) & $1.2\mathrm{e}{-1}$ (0) & $2.4\mathrm{e}{-1}$ (0) & $4.3\mathrm{e}{-1}$ (0) & $4.3\mathrm{e}{-1}$ (0) & $4.1\mathrm{e}{-1}$ (0) & $2.5\mathrm{e}{-1}$ (0) & $1.3\mathrm{e}{-1}$ (0) \\
 & 60 & $1.6\mathrm{e}{-1}$ (0) & $2.5\mathrm{e}{-3}$ (5) & $1.1\mathrm{e}{-1}$ (0) & $9.1\mathrm{e}{-2}$ (0) & $9.1\mathrm{e}{-2}$ (0) & $1.2\mathrm{e}{-1}$ (0) & $6.6\mathrm{e}{-2}$ (0) & $2.3\mathrm{e}{-11}$ (6) \\
 & 90 & $1.1\mathrm{e}{-1}$ (0) & $2.6\mathrm{e}{-11}$ (10) & $7.1\mathrm{e}{-2}$ (0) & $5.3\mathrm{e}{-2}$ (0) & $5.3\mathrm{e}{-2}$ (0) & $7.1\mathrm{e}{-2}$ (0) & $1.5\mathrm{e}{-2}$ (0) & $2.2\mathrm{e}{-11}$ (10) \\
\hline
EXP & 40 & $8.1\mathrm{e}{-1}$ & $2.9\mathrm{e}{-2}$ & $5.4\mathrm{e}{-1}$ & $3.0\mathrm{e}{-1}$ & $3.0\mathrm{e}{-1}$ & $3.0\mathrm{e}{-1}$ & $5.4\mathrm{e}{-1}$ & $3.0\mathrm{e}{-2}$ \\
 & 60 & $4.1\mathrm{e}{-1}$ & $1.6\mathrm{e}{-2}$ & $2.5\mathrm{e}{-1}$ & $6.1\mathrm{e}{-1}$ & $6.1\mathrm{e}{-1}$ & $5.7\mathrm{e}{-1}$ & $3.0\mathrm{e}{-1}$ & $2.1\mathrm{e}{-2}$ \\
 & 90 & $2.3\mathrm{e}{-1}$ & $2.6\mathrm{e}{-3}$ & $1.2\mathrm{e}{-1}$ & $3.0\mathrm{e}{-2}$ & $3.0\mathrm{e}{-2}$ & $5.8\mathrm{e}{-2}$ & $5.3\mathrm{e}{-2}$ & $1.7\mathrm{e}{-3}$ \\
\hline
\end{tabular}
\caption{Sparse polynomial approximation: median relative test error over 10 instances (for
ROS, number of exact recoveries in parentheses).}
\label{tabsparsepol}
\end{table}

In this problem the monomials are strongly correlated on the sampling region, and the
Telescopic Method started at $c = 0$ stops at B-stationary points that are frequently far
from the best approximations. The $\ell_1$ solutions approximate the functions well, but they do not satisfy the
cardinality constraints, except in the ROS instances in which the true polynomial is
recovered (5 instances with $m = 60$ and all the instances with $m = 90$). Started at the
$\ell_1$ solution, the Telescopic Method produces a polynomial that satisfies the cardinality
constraints with a median test error of the same order as the one of the $\ell_1$ solution:
slightly larger for ROS with $m = 40$ and for EXP with $m \in \{40, 60\}$, and smaller in the
other cases. In particular, for ROS with $m = 60$ the true polynomial was recovered in 6
instances (instead of 5) and the median test error decreased from $2.5 \times 10^{-3}$ to
$2.3 \times 10^{-11}$, and for EXP with $m = 90$ the median test error decreased from
$2.6 \times 10^{-3}$ to $1.7 \times 10^{-3}$. Started at $c = 0$, the Telescopic Method has smaller median
test errors than MN and HT except for ROS with $m = 40$ and for EXP with $m = 60$; started at
the HT point, it improves HT for ROS with $m \in \{60, 90\}$ and for EXP with $m = 90$. In this problem
the scaled coefficients are small, so that the choice $\delta = 0.1$ has almost no effect;
the look-ahead strategy reduced the median test error in three of the four cases with
$m \geq 60$. The median CPU time of the Telescopic Method
was about 0.04 seconds per run. In all the runs of the Telescopic Method (all the
versions, HT+T and L1+T), the method stopped at an approximately B-stationary point, with the
B-stationarity measure~(\ref{bstatmeas}) at most $10^{-8}$.

\paragraph{Local minimizers and interpolation.}
For ROS, the true polynomial satisfies the cardinality constraints and interpolates the data,
so that it is a global minimizer of the PCC problem, with objective value zero. Started at
$c = 0$ or at the HT point, the Telescopic Method never reached it: it stopped at B-stationary
points with relative training errors around $10^{-2}$. Started at the $\ell_1$ solution, it
reached the global minimizer in 6 and 10 instances for $m = 60$ and $m = 90$, respectively. As
in separable regression, the difficulty lies in the existence of non-global B-stationary
points, and the initial point is decisive.

The cardinality constraints act as a regularization only when the sparse model cannot
interpolate the data. For EXP with $m = 40$, the number of coefficients allowed by the
cardinality constraints ($7 + 11 + 16 + 21 = 55$) exceeds the number of sample points, and the
Telescopic Method started at $c = 0$ and at the HT point found interpolating polynomials (relative
training errors below $10^{-8}$) in all the instances. All of them are global minimizers of the
PCC problem, but their median test errors (0.30 and 0.54) are much larger than the one of
L1+T (0.030). In this case, the objective function does not discriminate between supports;
choosing, among seven runs (from $c = 0$, from the HT and $\ell_1$ points, and from four random
faces), the solution with smallest objective value gives a median test error of 0.54, and the
quality of the result is determined by the initial point. When the sparse model cannot
interpolate the data (ROS for all $m$ and EXP for $m \geq 60$), the training error is a good
guide: among the solutions of T, HT+T, and L1+T, the one with smallest training error had the
smallest test error in 45 of the 50 instances, and the best of the seven runs was the L1+T
solution in 49 of them.

\subsection{Portfolio selection with group cardinality constraints} \label{numerportfolio}

We consider $n$ assets divided into $N$ groups (e.g., economic sectors) with $n_\nu$ assets in
group $\nu$. The current price of asset $j$ is $\bar\theta_j$ and, under scenario
$i \in \{1, \ldots, m\}$, its price at the end of the investment period is $\theta_{ij}$. Let
$r_{ij} = \theta_{ij}/\bar\theta_j - 1$ be the corresponding return. An investor wants to
invest an amount $M$; the decision variable $x_j$ is the amount invested in asset $j$.
Under scenario $i$ the profit is $R_i(x) = \sum_{j=1}^n r_{ij} x_j$, with mean and standard
deviation
\[
\mu(x) = \frac{1}{m} \sum_{i=1}^m R_i(x), \qquad
\sigma(x) = \left( \frac{1}{m} \sum_{i=1}^m \left[ R_i(x) - \mu(x) \right]^2 \right)^{1/2}.
\]
A rigid budget constraint $\sum_j x_j = M$ (or $\leq M$) is not realistic, because the
investor may borrow money if this is worthwhile, or keep part of the money in a risk-free
investment. We model the cost of borrowing and the gain of the risk-free investment by
quadratic terms and, as in the classical model of Markowitz \cite{markowitz}, we measure the
risk by the variance $\sigma(x)^2$. This leads to
\begin{equation} \label{portf}
\Minimize \; -\mu(x) + \omega \, \sigma(x)^2
  + \rho \Big[ \Big( \sum_{j=1}^n x_j - M \Big)_+ \Big]^2
  - r \Big[ \Big( M - \sum_{j=1}^n x_j \Big)_+ \Big]^2 ,
\end{equation}
where $(s)_+ = \max\{0, s\}$, $\omega \geq 0$ reflects risk aversion, and $\rho > r \geq 0$.
The squares in the last two terms make them continuously differentiable with Lipschitz
continuous derivatives. Short selling is excluded by $0 \leq x_j \leq u_j$; since the
$u_j$ are finite, $\B$ is bounded and Assumptions~\ref{abox}--\ref{alevel} hold. Finally,
regulations (or a conservative investor) impose that at most $p_\nu$ assets of group $\nu$ are
used, which gives a PCC problem. In \cite{bmaus}, a cardinality constraint on the whole
portfolio was handled by a sequence of smooth approximations of the step function; here, the
cardinality constraints of each group are handled exactly. (Mean-variance models with a
single cardinality constraint were addressed, for example, in \cite{bienstock,cmbs00,kkr23}.)

\paragraph{Why the variance.}
A model with the standard deviation $\sigma(x)$, instead of $\sigma(x)^2$, in~(\ref{portf})
would be natural, since $\sigma(x)$ is measured in the same units as $\mu(x)$. However,
$\sigma$ is not differentiable at the points in which all the returns $R_i(x)$ coincide, in
particular at $x = 0$ (the portfolio that keeps the whole capital in the risk-free investment),
where it grows linearly along every direction $d$ with $\sigma(d) > 0$. As a consequence, for large values of $\omega$
the point $x = 0$ is a local minimizer of the (nonsmooth) objective function on faces of the
feasible set, and, even when it is not, methods designed for smooth problems may stop there.
In both cases B-stationarity cannot be verified at $x = 0$, because the objective function is
not differentiable there. This is exactly what we observed in preliminary experiments with
$\sigma(x)$ (see the end of this section). With the
variance, the objective function is smooth; moreover, its gradient at $x = 0$ is
$-\bar r + 2 r M e$, where $\bar r_j = (1/m) \sum_i r_{ij}$, so
that $x = 0$ is not B-stationary (and, hence, not a local minimizer) whenever the expected
return of some asset exceeds $2 r M$.

\paragraph{Choice of the weight $\omega$.}
The two risk measures have different scales: $\sigma(tx) = t \sigma(x)$, whereas
$\sigma(tx)^2 = t^2 \sigma(x)^2$. In a model with $\omega_1 \sigma(x)$, one additional unit of
standard deviation costs $\omega_1$ units of expected return, whereas in~(\ref{portf}) it costs
$d(\omega \sigma^2)/d\sigma = 2 \omega \sigma$, which depends on the level of risk. We choose
$\omega$ so that both models attribute the same marginal price to risk at a reference level
$\sigma_{\mathrm{ref}}$, that is,
\begin{equation} \label{omegacal}
\omega = \frac{\omega_1}{2 \sigma_{\mathrm{ref}}},
\end{equation}
where $\sigma_{\mathrm{ref}}$ is the standard deviation of a typical portfolio (in the
experiments, the equal-weight portfolio). This choice has a classical interpretation: when the
returns are normally distributed, $\mu(x) - (\lambda/2) \sigma(x)^2$ is the certainty equivalent
of the return of an investor with constant absolute risk aversion $\lambda$; since the capital
is normalized ($M = 1$), $\lambda = 2 \omega$ is the relative risk aversion of the investor, for
which values not larger than 10 are usually considered reasonable \cite{mehraprescott}.

\paragraph{Data and settings.}
We used the daily adjusted closing prices of the
constituents of the S\&P~500 index in the 503 trading days from October 16, 2024, to
October 19, 2026 (Yahoo Finance), and their GICS sectors (Wikipedia). The data were
downloaded on October 20, 2026, with the Python package \texttt{yfinance}. Since the current
constituents of the index were used, assets that left the index during the period are not
considered (survivorship bias), which contributes to the large expected returns discussed
below. The 497 assets with
complete data form $N = 11$ groups (the sectors), with 21 to 81 assets each. We generated
$m = 10{,}000$ scenarios (12 of them catastrophic), as described below,
and used $M = 1$, $r = 0.04$, and $p_\nu = p \in \{1, 2, 3\}$ for all $\nu$. Two settings
were considered: (A) at most 20\% of the capital in each asset ($u_j = 0.2$) and moderate
borrowing cost ($\rho = 1$); and (B) no bound on the fraction of the capital invested in an
asset ($u_j = 1$) and expensive borrowing ($\rho = 10$). The standard deviation of the
equal-weight portfolio is $\sigma_{\mathrm{ref}} = 0.2076$ and, according
to~(\ref{omegacal}), we used $\omega = \omega_1 / (2 \sigma_{\mathrm{ref}})$ with
$\omega_1 \in \{0.5, 1, 2\}$, which gives $\omega \in \{1.20, 2.41, 4.82\}$ and relative risk
aversions $\lambda = 2 \omega \in \{2.4, 4.8, 9.6\}$.

The scenarios are generated as in \cite{bmaus}: from the historical daily adjusted closing
prices $d_{t,j}$ of the assets ($t$ is the day and $j$ the asset), the daily growth factors
$d_{t+1,j}/d_{t,j}$ are computed; each
scenario is obtained by multiplying 250 growth factors (one year of trading days) taken from
randomly chosen days, so that $\theta_{ij}/\bar\theta_j$ is the product of these factors; and,
with probability $0.001$, a scenario is multiplied by a random number in $[0.1, 0.4]$ in order
to simulate catastrophic events.

The Telescopic Method (T and its variants) was started at the portfolio that invests
$M/(N p)$ (or $u_j$, if smaller) in each of the $p$ assets with largest expected return
$\bar r_j$ of each sector. As for separable regression, we also considered the problem without
cardinality constraints (REL), solved by GENCAN from the equal-weight portfolio, the hard
thresholding strategy HT (the $p$ largest components of the REL solution in each sector are
kept and the objective function is minimized by GENCAN on the corresponding face), and the
Telescopic Method started at the HT point (HT+T). The results are in
Table~\ref{tabportf}. The version with $\delta = 0.1$, which is not shown, found the same
solutions as T in 16 of the 18 instances and slightly worse ones in the other two (in setting
B, $f = -0.5371$ instead of $-0.5384$ for $\lambda = 4.8$ and $p = 2$, and $f = -0.3742$
instead of $-0.3757$ for $\lambda = 9.6$ and $p = 3$).

\begin{table}[ht!]
\centering
\scriptsize
\setlength{\tabcolsep}{2.5pt}
\begin{tabular}{cccccccc}
\hline
Set. & $\omega$ ($\lambda$) & REL & $p$ & HT & T & T (no l.-a.) & HT+T \\
\hline
A & 1.20 (2.4) & $-0.7947$ (19, 9) & 1 & $-0.6753$ (8) & $-0.6836$ (10) & $-0.6836$ (10) & $\mathbf{-0.6867}$ (10) \\
 &  &  & 2 & $-0.7609$ (14) & $-0.7578$ (15) & $\mathbf{-0.7644}$ (14) & $-0.7609$ (14) \\
 &  &  & 3 & $\mathbf{-0.7787}$ (15) & $-0.7776$ (14) & $-0.7659$ (15) & $\mathbf{-0.7787}$ (15) \\
\cline{2-8}
 & 2.41 (4.8) & $-0.5742$ (25, 9) & 1 & $-0.5033$ (9) & $\mathbf{-0.5142}$ (10) & $\mathbf{-0.5142}$ (10) & $-0.5107$ (10) \\
 &  &  & 2 & $-0.5480$ (15) & $-0.5526$ (17) & $\mathbf{-0.5534}$ (17) & $-0.5529$ (17) \\
 &  &  & 3 & $-0.5592$ (21) & $\mathbf{-0.5629}$ (20) & $\mathbf{-0.5629}$ (20) & $-0.5617$ (21) \\
\cline{2-8}
 & 4.82 (9.6) & $-0.3926$ (32, 10) & 1 & $-0.3221$ (10) & $\mathbf{-0.3422}$ (10) & $-0.3412$ (10) & $\mathbf{-0.3422}$ (10) \\
 &  &  & 2 & $-0.3660$ (17) & $-0.3746$ (18) & $-0.3746$ (18) & $\mathbf{-0.3769}$ (16) \\
 &  &  & 3 & $\mathbf{-0.3832}$ (20) & $-0.3821$ (21) & $-0.3819$ (21) & $\mathbf{-0.3832}$ (20) \\
\hline
B & 1.20 (2.4) & $-0.7640$ (17, 8) & 1 & $-0.6952$ (9) & $-0.6861$ (10) & $-0.6847$ (10) & $\mathbf{-0.7061}$ (9) \\
 &  &  & 2 & $\mathbf{-0.7324}$ (12) & $-0.7152$ (13) & $\mathbf{-0.7324}$ (12) & $\mathbf{-0.7324}$ (12) \\
 &  &  & 3 & $\mathbf{-0.7507}$ (11) & $-0.7464$ (13) & $-0.7464$ (13) & $\mathbf{-0.7507}$ (11) \\
\cline{2-8}
 & 2.41 (4.8) & $-0.5564$ (22, 9) & 1 & $-0.5025$ (10) & $-0.4836$ (10) & $-0.4801$ (10) & $\mathbf{-0.5032}$ (10) \\
 &  &  & 2 & $-0.5364$ (15) & $\mathbf{-0.5384}$ (16) & $-0.5371$ (16) & $-0.5370$ (16) \\
 &  &  & 3 & $-0.5414$ (18) & $\mathbf{-0.5459}$ (18) & $\mathbf{-0.5459}$ (18) & $-0.5435$ (18) \\
\cline{2-8}
 & 4.82 (9.6) & $-0.3880$ (29, 10) & 1 & $-0.3202$ (10) & $-0.3270$ (10) & $-0.3411$ (10) & $\mathbf{-0.3413}$ (10) \\
 &  &  & 2 & $-0.3626$ (17) & $-0.3644$ (17) & $-0.3644$ (17) & $\mathbf{-0.3692}$ (16) \\
 &  &  & 3 & $\mathbf{-0.3789}$ (20) & $-0.3757$ (22) & $-0.3757$ (22) & $\mathbf{-0.3789}$ (20) \\
\hline
\end{tabular}
\caption{Portfolio selection with 497 assets of the S\&P~500 in 11 sectors: objective
function value and, in parentheses, number of assets in the portfolio. For REL (no
cardinality constraints) the second number in parentheses is the largest number of assets of
a sector in the portfolio. In each row, the best value among the methods that satisfy the
cardinality constraints is in bold.}
\label{tabportf}
\end{table}

All the runs of the Telescopic Method stopped with the practical stopping criterion
satisfied, with B-stationarity measure below $10^{-14}$, after at most 16 swaps, and the CPU
time of every run (including, for HT and HT+T, the solution of REL) was smaller than one
second. The solutions of REL use between 17 and 32 assets, with up to 10 assets in a sector,
so that the cardinality constraints are far from being satisfied. Started at the
portfolio described above, the Telescopic Method found a smaller objective function value
than HT in 10 of the 18 instances. Started at the HT point, it improved the HT solution in 12
instances; in the other 6 instances, it ended with the same objective function value as
HT. The best of the T and HT+T solutions was at least as good as
the HT solution in all the instances. This illustrates the main practical difference between
the Telescopic Method and heuristics such as HT: the latter produce a point that is a
minimizer on a face chosen once and for all, whereas the former keeps exchanging assets
within the sectors while this is profitable, and stops only at a B-stationary point. The
look-ahead strategy improved the objective function value in 6 instances and worsened it in
4.

\begin{figure}[ht!]
\centering
\includegraphics[width=\textwidth]{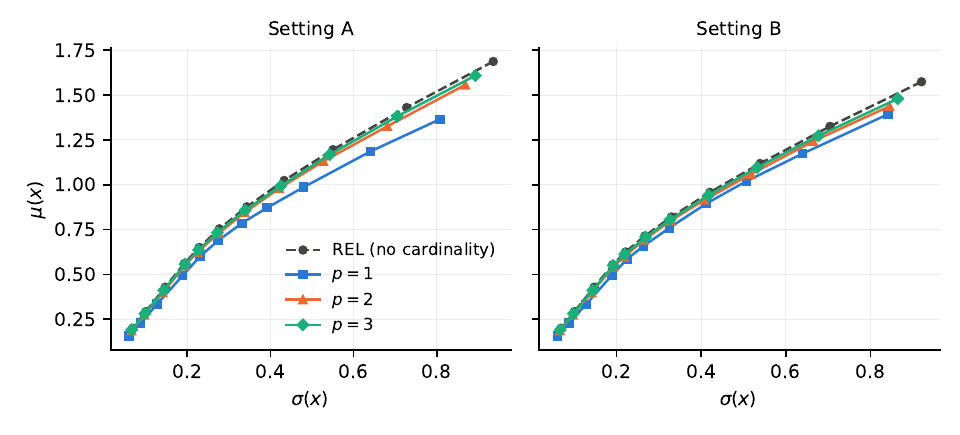}
\caption{Expected return versus standard deviation of the portfolios obtained with
$\omega = \omega_1 / (2 \sigma_{\mathrm{ref}})$, $\omega_1 \in \{0.25, 0.35, 0.5, 0.7, 1, 1.4, 2,
2.8, 4, 5.6, 8\}$, without (REL) and with ($p = 1, 2, 3$) the sector cardinality constraints.
The risk aversion decreases from left to right.}
\label{figfront}
\end{figure}

\paragraph{The price of the cardinality constraints.}
Table~\ref{tabportfchar} describes the portfolios for $\lambda = 4.8$, and
Figure~\ref{figfront} shows the expected return $\mu(x)$ and the standard deviation
$\sigma(x)$ of the portfolios obtained for a grid of eleven values of $\omega_1$ between
$0.25$ and $8$ (that is, $\lambda$ between $1.2$ and $38.5$); for $p = 1, 2, 3$ we report the
best of the T and HT+T solutions. Since the objective function is the negative of the
risk-adjusted return, its relative increase measures what the investor loses by imposing
the cardinality constraints. In the 18 instances of Table~\ref{tabportf}, this loss was
between 7.6\% and 13.6\% for $p = 1$, between 3.2\% and 4.8\% for $p = 2$, and between
1.7\% and 2.4\% for $p = 3$, whereas the number of assets in the portfolio decreased from 17--32
(REL) to 9--10, 12--17, and 11--20, respectively. Thus, allowing two or three assets per
sector gives portfolios that are almost as good as the unconstrained ones, as
Figure~\ref{figfront} confirms: the curves for $p = 2$ and $p = 3$ are very close to the one
of REL, and only the curve for $p = 1$ is visibly below it, mainly for small risk aversion. For
$\lambda = 4.8$ and $\lambda = 9.6$, the standard deviation is almost unaffected by the
cardinality constraints (for $\lambda = 4.8$ it varies between $0.325$ and $0.343$ in
Table~\ref{tabportfchar}), and the loss comes from the expected return; for $\lambda = 2.4$,
both the expected return and the standard deviation decrease when $p = 1$ is imposed.

The constraints also change the composition of the portfolios. Without them, the portfolio
is concentrated in Information Technology (8 to 10 assets) and Health Care (3 to 7 assets),
uses no asset of the Utilities sector and, except for $\lambda = 9.6$, no asset of the
Financials sector, and covers 7 to 10 of the 11 sectors. With $p = 1$, the portfolio covers 10
sectors in all the instances except one. In this sense, the group cardinality constraints act
as diversification constraints. On the other hand, the remaining assets receive larger
weights: in setting B, which has no bound per asset, the largest weight increases from 0.32
to 0.37 for $\lambda = 2.4$, from 0.20 to 0.27 for $\lambda = 4.8$, and from 0.13 to 0.23 for
$\lambda = 9.6$ when $p = 1$ is imposed.

\paragraph{Settings A and B.}
In setting A, in which borrowing is cheap but at most 20\% of the capital can be invested in
an asset, the investor borrows between 4\% and 19\% of the capital and the bound $u_j = 0.2$ is
active for several assets. In setting B, in which borrowing is expensive but there is no bound
per asset, the borrowing is at most 2\% and the portfolios are more concentrated (the largest
weight reaches 0.37 for $\lambda = 2.4$ and $p = 1$, and 0.54 in Figure~\ref{figfront} for
$\lambda = 1.2$). Cheap leverage pays off: the objective function values are slightly smaller
in setting A in all the instances of Table~\ref{tabportf} except one. For large risk aversion ($\lambda \geq 19.3$ in Figure~\ref{figfront}) the investor keeps
part of the capital in the risk-free investment, so that neither the bounds nor the cost of
borrowing are active, and the two settings give the same portfolios, except for $p = 1$ and
$\lambda = 19.3$, in which the method stopped at different B-stationary points in the two
settings ($f = -0.1833$ in setting A and $-0.1790$ in setting B). In this regime, the
cardinality constraints also reduce the invested fraction (for $\lambda = 19.3$, from 83\%
without them to 60\% and 68\% with $p = 1$ in settings A and B, respectively), since fewer
assets give less diversification.

\begin{table}[ht!]
\centering
\scriptsize
\setlength{\tabcolsep}{3pt}
\begin{tabular}{cccccccccc}
\hline
Set. & & $f$ & loss (\%) & $\mu$ & $\sigma$ & $\sum_j x_j$ & assets & sectors & largest weight \\
\hline
A & REL & $-0.5742$ & -- & 0.877 & 0.343 & 1.135 & 25 & 8 & 0.200 (VLO) \\
 & $p = 1$ & $-0.5142$ & 10.5 & 0.785 & 0.331 & 1.081 & 10 & 10 & 0.200 (WELL) \\
 & $p = 2$ & $-0.5534$ & 3.6 & 0.845 & 0.338 & 1.126 & 17 & 10 & 0.200 (VLO) \\
 & $p = 3$ & $-0.5629$ & 2.0 & 0.858 & 0.339 & 1.133 & 20 & 10 & 0.200 (VLO) \\
\hline
B & REL & $-0.5564$ & -- & 0.821 & 0.330 & 1.015 & 22 & 8 & 0.196 (VLO) \\
 & $p = 1$ & $-0.5032$ & 9.6 & 0.760 & 0.325 & 1.012 & 10 & 10 & 0.270 (CAH) \\
 & $p = 2$ & $-0.5384$ & 3.2 & 0.800 & 0.328 & 1.014 & 16 & 10 & 0.235 (CAH) \\
 & $p = 3$ & $-0.5459$ & 1.9 & 0.805 & 0.327 & 1.014 & 18 & 9 & 0.210 (VLO) \\
\hline
\end{tabular}
\caption{Portfolios for $\lambda = 4.8$ ($\omega = 2.41$): objective function value, its
relative increase with respect to REL, expected return, standard deviation, invested
fraction of the capital, number of assets, number of sectors, and largest weight. For
$p = 1, 2, 3$, the best solution found by the methods that satisfy the cardinality constraints
is described.}
\label{tabportfchar}
\end{table}

\paragraph{In-sample and out-of-sample returns.}
The expected returns in Table~\ref{tabportfchar} (around 80\% a year) are much larger than
what an investor may reasonably expect. They reflect the exceptional performance of some
assets in the two years of data (for example, the five largest values of $\bar r_j$ are between
2.2 and 6.9), which the scenarios, built by resampling the observed daily growth factors,
reproduce. Hence, these values should be read as in-sample quantities. To get an idea of the
out-of-sample behavior, we repeated the experiments building the scenarios with the first
year of data only (October 16, 2024, to October 17, 2025), with $\sigma_{\mathrm{ref}}$
recomputed from these scenarios, and computed the realized return
$\sum_j x_j (d_{T,j}/d_{h,j} - 1)$ of each portfolio over the second year, where $h$ and $T$
are the first and the last day of the second year. Table~\ref{tabportfoos} shows that the
realized returns (between 7\% and 33\%; 18\% for the equal-weight portfolio) are much smaller
than the expected ones (between 72\% and 163\%). With $p = 2$ and $p = 3$, the realized
returns are close to the ones of REL (the largest difference is 0.11, for setting A and
$\omega_1 = 0.5$), whereas with $p = 1$ they are smaller in four of the six cases, essentially
equal in one, and larger in one. Since this is a single realization of one year, these numbers
are only indicative; they suggest, however, that the in-sample losses due to the cardinality
constraints with $p \geq 2$ are small compared with the uncertainty of the estimates, and that
these constraints did not deteriorate the portfolios out of sample.

\begin{table}[ht!]
\centering
\scriptsize
\begin{tabular}{cccccc}
\hline
Set. & $\omega_1$ & REL & $p = 1$ & $p = 2$ & $p = 3$ \\
\hline
A & 0.5 & 1.63 / 0.30 & 1.51 / 0.07 & 1.58 / 0.19 & 1.63 / 0.33 \\ 
A & 1.0 & 1.12 / 0.24 & 1.06 / 0.22 & 1.11 / 0.24 & 1.12 / 0.25 \\ 
A & 2.0 & 0.81 / 0.21 & 0.72 / 0.26 & 0.77 / 0.23 & 0.80 / 0.21 \\ 
B & 0.5 & 1.54 / 0.29 & 1.49 / 0.10 & 1.53 / 0.28 & 1.54 / 0.29 \\ 
B & 1.0 & 1.05 / 0.26 & 1.02 / 0.15 & 1.05 / 0.27 & 1.06 / 0.28 \\ 
B & 2.0 & 0.77 / 0.21 & 0.78 / 0.22 & 0.79 / 0.22 & 0.77 / 0.22 \\ 
\hline
\end{tabular}
\caption{Expected return (scenarios built with the first year of data) / realized return over
the second year.}
\label{tabportfoos}
\end{table}

\paragraph{Preliminary experiments with $\sigma(x)$.}
Before adopting the variance, we ran the same experiments with $\omega_1 \sigma(x)$ instead of
$\omega \sigma(x)^2$ in~(\ref{portf}), with $\omega_1 \in \{0.5, 1, 2\}$. Where $\sigma$ is not
differentiable, the implementation uses the subgradient $0$ of $\sigma$. With $\omega_1 = 2$,
the Telescopic Method started at the portfolio described above (T) converged to $x = 0$, whose
objective function value is $-r M^2 = -0.04$, in 2 of the 6 instances (setting A with
$p = 2$ and $p = 3$), and so did the versions with $\delta = 0.1$ and without look-ahead. In
all these runs the method stopped because no swap decreased $f$ although the swap-candidates
test was not satisfied. Indeed, the ``gradient''
used at $x = 0$, namely $-\bar r + 2 r M e$, indicates that adding any asset with
$\bar r_j > 2 r M$ decreases $f$, but the directional derivative of the actual objective
function along $d \geq 0$ is $(-\bar r + 2 r M e)^T d + \omega_1 (d^T C d)^{1/2}$, which, for
$\omega_1 = 2$, is positive along every coordinate direction $d = e_j$. Thus, no single asset
can be added with profit. However, in both cases there are directions of the final face
(diversified portfolios) along which $f$ decreases, but GENCAN, which assumes
differentiability, stopped at $x = 0$. In any case, B-stationarity cannot be verified at
$x = 0$. In the same instances, HT and HT+T found portfolios with 17 to 20 assets, investing
(slightly) more than the capital, with objective function values between $-0.18$ and $-0.16$,
much smaller than $-0.04$. With the variance, no run of any method ended at $x = 0$; in fact, in all the
runs the whole capital was invested.

\section{Conclusions} \label{concl}

We introduced the Telescopic Method for the minimization of a smooth function subject to
bounds and partitioned cardinality constraints, a class that contains bound-constrained
minimization, cardinality-constrained minimization with bounds, and bound-constrained
mathematical programs with complementarity constraints. The method alternates iterations
on the current face, of any bound-constrained method that satisfies the hypotheses of
Section~\ref{facemin}, with one-dimensional escaping tests that exchange a
variable of the face by a variable outside it within the same group. We proved that every
limit point of the generated sequence is B-stationary. This result requires neither
constraint qualifications nor global or exact solutions of subproblems, and the method does
not solve linear programs, or programs with complementarity constraints, at any iteration.
In contrast, the convergence results of the methods based on the continuous reformulation of
cardinality constraints concern stationary points of the reformulated problem, which need not
be B-stationary for the original one.

The method was implemented in Fortran with GENCAN as the method for the minimization on the
faces, using exact second derivatives. In the 60 bound-constrained MPCCs of the test set
of Kirches, Larson, Leyffer and Manns, it found approximately B-stationary points in all the
problems. In the
applications, which involve groups with more than two variables and more than one non-null
variable per group, the method was particularly useful as a refinement of the points given
by hard thresholding or by $\ell_1$ minimization. In separable regression, started at the hard thresholding point,
the method recovered the true support in instances in which thresholding failed. In sparse
polynomial approximation, started at the $\ell_1$ solution, it produced polynomials that
satisfy the cardinality constraints with test errors of the same order as the ones of the
$\ell_1$ solutions. In portfolio selection with up to $p$ assets per sector, with real data
of 497 assets of the S\&P~500, the best of the solutions obtained from the two initial points
was at least as good as the hard thresholding solution in all the instances; moreover,
allowing two or three assets per sector cost only 2\% to 5\% of the risk-adjusted return of the
unconstrained portfolios, while spreading the capital over more sectors. The experiments
also showed that letting every variable of the current face be a candidate to leave it
($\delta = +\infty$) is much better than restricting the candidates to the variables with
small absolute value in the MPCCs and in separable regression, and makes little or no
difference in the other applications, and that the look-ahead strategy, with GENCAN run to
convergence on the candidate faces, improves the results more often than it worsens them, at
a moderate cost (in the MPCCs, it roughly doubled the CPU time). Runs from several initial points showed that the main
difficulty in the applications is not the convergence to B-stationary points, which was always
achieved, but the existence of B-stationary points that are not global minimizers; good
initial points, such as the ones given by hard thresholding or by $\ell_1$ minimization, were
decisive. In the MPCCs, keeping the best of 21 runs matched or improved the best
published values in 55 of the 60 problems.

The portfolio experiments illustrated a limitation that is inherent to methods based on
derivatives: with the standard deviation as risk measure, the objective function
is not differentiable at $x = 0$, and, in preliminary experiments, the method stopped there in
some instances, a point at which B-stationarity cannot be verified. Using the variance, with a weight calibrated to
preserve the marginal price of risk, removed this difficulty.

Problems with additional equality and inequality constraints, such as mathematical programs
with complementarity constraints with general constraints, can be handled by an augmented
Lagrangian method in which the bounds and the partitioned cardinality constraints are kept as
lower-level constraints, so that the subproblems are solved by the Telescopic Method. This
is the subject of a work in progress \cite{bmal}. In that framework, B-stationarity of the
limit points can be guaranteed, without constraint qualifications, if the penalty parameters
remain bounded, but it may be lost otherwise, because a swap passes through infeasible points
whose penalty grows with the penalty parameter. We are studying swaps that preserve the
values of the additional constraints, with which the swap tests do not depend on the penalty
parameter.

Other extensions are possible. Since the escaping tests are independent, they can be
performed in parallel. Finally, Assumption~\ref{abox} requires that zero belongs to the bounds of the
variables of the constrained blocks; the treatment of semicontinuous variables, for which
non-null values must lie in an interval that does not contain zero, would require a
different escaping procedure.

\section*{Declaration of generative AI use}
During the preparation of this work, the authors used Claude (Anthropic) to assist in drafting the
manuscript, checking the proofs, implementing the codes and running the numerical
experiments. The authors reviewed and edited all the content and take full responsibility for the
content of the publication.



\begin{thebibliography}{99}

\bibitem{afss19} R. Andreani, N. S. Fazzio, M. L. Schuverdt, and L. D. Secchin, A sequential
  optimality condition related to the quasi-normality constraint qualification and its
  algorithmic consequences, {\it SIAM Journal on Optimization} 29, pp.~743--766, 2019.

\bibitem{ahm11} R. Andreani, G. Haeser, and J. M. Mart\'{\i}nez, On sequential optimality
  conditions for smooth constrained optimization, {\it Optimization} 60, pp.~627--641, 2011.

\bibitem{ahss19} R. Andreani, G. Haeser, L. D. Secchin, and P. J. S. Silva, New sequential
  optimality conditions for mathematical programs with complementarity constraints and
  algorithmic consequences, {\it SIAM Journal on Optimization} 29, pp.~3201--3230,
  2019.

\bibitem{ams10} R. Andreani, J. M. Mart\'{\i}nez, and B. F. Svaiter, A new sequential optimality
  condition for constrained optimization and algorithmic consequences, {\it SIAM Journal on
  Optimization} 20, pp.~3533--3554, 2010.

\bibitem{ass18} R. Andreani, L. D. Secchin, and P. J. S. Silva, Convergence properties of a second
  order augmented Lagrangian method for mathematical programs with complementarity constraints,
  {\it SIAM Journal on Optimization} 28, pp.~2574--2600, 2018.

\bibitem{abm} M. Andretta, E. G. Birgin, and J. M. Mart\'{\i}nez, Practical active-set Euclidean
  trust-region method with spectral projected gradients for bound-constrained minimization,
  {\it Optimization} 54, pp.~305--325, 2005.

\bibitem{atw07} M. Anitescu, P. Tseng, and S. J. Wright, Elastic-mode algorithms for mathematical
  programs with equilibrium constraints: global convergence and stationarity properties,
  {\it Mathematical Programming} 110, pp.~337--371, 2007.

\bibitem{bsv2012} A. S. Bandeira, K. Scheinberg, and L. N. Vicente, Computation of sparse low degree
  interpolating polynomials and their application to derivative-free optimization,
  {\it Mathematical Programming} 134, pp.~223--257, 2012.

\bibitem{be13} A. Beck and Y. C. Eldar, Sparsity constrained nonlinear optimization: optimality
  conditions and algorithms, {\it SIAM Journal on Optimization} 23, pp.~1480--1509, 2013.

\bibitem{bh16} A. Beck and N. Hallak, On the minimization over sparse symmetric sets: projections,
  optimality conditions, and algorithms, {\it Mathematics of Operations Research} 41,
  pp.~196--223, 2016.

\bibitem{bienstock} D. Bienstock, Computational study of a family of mixed-integer quadratic
  programming problems, {\it Mathematical Programming} 74, pp.~121--140, 1996.

\bibitem{gencan} E. G. Birgin and J. M. Mart\'{\i}nez, Large-scale active-set box-constrained
  optimization method with spectral projected gradients, {\it Computational Optimization and
  Applications} 23, pp.~101--125, 2002.

\bibitem{bmbook} E. G. Birgin and J. M. Mart\'{\i}nez, {\it Practical Augmented Lagrangian Methods
  for Constrained Optimization}, SIAM, Philadelphia, 2014.

\bibitem{bmaus} E. G. Birgin and J. M. Mart\'{\i}nez, On the application of an Augmented Lagrangian
  algorithm to some portfolio problems, {\it EURO Journal on Computational Optimization} 4,
  pp.~79--92, 2016.

\bibitem{bm2020} E. G. Birgin and J. M. Mart\'{\i}nez, Complexity and performance of an Augmented
  Lagrangian algorithm, {\it Optimization Methods and Software} 35, pp.~885--920, 2020.

\bibitem{bmal} E. G. Birgin and J. M. Mart\'{\i}nez, Augmented Lagrangian methods with
  telescopic subproblems for MPCCs and cardinality constraints with general constraints,
  work in progress.

\bibitem{spg} E. G. Birgin, J. M. Mart\'{\i}nez, and M. Raydan, Nonmonotone spectral projected
  gradient methods on convex sets, {\it SIAM Journal on Optimization} 10, pp.~1196--1211, 2000.

\bibitem{bd2009} T. Blumensath and M. E. Davies, Iterative hard thresholding for compressed sensing,
  {\it Applied and Computational Harmonic Analysis} 27, pp.~265--274, 2009.

\bibitem{bbcs18} M. Branda, M. Bucher, M. \v{C}ervinka, and A. Schwartz, Convergence of a
  Scholtes-type regularization method for cardinality-constrained optimization problems with an
  application in sparse robust portfolio optimization, {\it Computational Optimization and
  Applications} 70, pp.~503--530, 2018.

\bibitem{bks} O. P. Burdakov, C. Kanzow, and A. Schwartz, Mathematical programs with cardinality
  constraints: reformulation by complementarity-type conditions and a regularization method,
  {\it SIAM Journal on Optimization} 26, pp.~397--425, 2016.

\bibitem{cks16} M. \v{C}ervinka, C. Kanzow, and A. Schwartz, Constraint qualifications and
  optimality conditions for optimization problems with cardinality constraints,
  {\it Mathematical Programming} 160, pp.~353--377, 2016.

\bibitem{cmbs00} T.-J. Chang, N. Meade, J. E. Beasley, and Y. M. Sharaiha, Heuristics for
  cardinality constrained portfolio optimisation, {\it Computers \& Operations Research} 27,
  pp.~1271--1302, 2000.

\bibitem{fl04} R. Fletcher and S. Leyffer, Solving mathematical programs with complementarity
  constraints as nonlinear programs, {\it Optimization Methods and Software} 19, pp.~15--40,
  2004.

\bibitem{ft02} M. Fukushima and P. Tseng, An implementable active-set algorithm for computing a
  B-stationary point of a mathematical program with linear complementarity constraints,
  {\it SIAM Journal on Optimization} 12, pp.~724--739, 2002; Erratum, {\it SIAM Journal on
  Optimization} 17, pp.~1253--1257, 2007.

\bibitem{gkr25} R. Garmanjani, E. H. M. Krulikovski, and A. Ramos, On stationarity conditions and
  constraint qualifications for multiobjective optimization problems with cardinality constraints,
  {\it Applied Mathematics \& Optimization} 91, article 22, 2025.

\bibitem{hm20} H. Hazimeh and R. Mazumder, Fast best subset selection: coordinate descent and local
  combinatorial optimization algorithms, {\it Operations Research} 68, pp.~1517--1537, 2020.

\bibitem{hlzz26} Z. Hu, X. Liu, S. Zeng, and J. Zhang, Randomized branch methods with inexact
  subproblems for Bouligand stationarity in linear and quadratic programs with complementarity
  constraints, arXiv:2609.29888, 2026.

\bibitem{hh2018} Q. Huangfu and J. A. J. Hall, Parallelizing the dual revised simplex method,
  {\it Mathematical Programming Computation} 10, pp.~119--142, 2018.

\bibitem{isu12} A. F. Izmailov, M. V. Solodov, and E. I. Uskov, Global convergence of augmented
  Lagrangian methods applied to optimization problems with degenerate constraints, including
  problems with complementarity constraints, {\it SIAM Journal on Optimization} 22,
  pp.~1579--1606, 2012.

\bibitem{jkmw} X. Jia, C. Kanzow, P. Mehlitz, and G. Wachsmuth, An augmented Lagrangian method
  for optimization problems with structured geometric constraints, {\it Mathematical Programming}
  199, pp.~1365--1415, 2023.

\bibitem{jsrf07} J. J. J\'udice, H. D. Sherali, I. M. Ribeiro, and A. M. Faustino, Complementarity
  active-set algorithm for mathematical programming problems with equilibrium constraints,
  {\it Journal of Optimization Theory and Applications} 134, pp.~467--481, 2007.

\bibitem{krsal21} C. Kanzow, A. B. Raharja, and A. Schwartz, An augmented Lagrangian method for
  cardinality-constrained optimization problems, {\it Journal of Optimization Theory and
  Applications} 189, pp.~793--813, 2021.

\bibitem{kllm} C. Kirches, J. Larson, S. Leyffer, and P. Manns, Sequential linearization method for
  bound-constrained mathematical programs with complementarity constraints, {\it SIAM Journal on
  Optimization} 32, pp.~75--99, 2022.

\bibitem{kllmtr} C. Kirches, J. Larson, S. Leyffer, and P. Manns, Sequential linearization method
  for bound-constrained mathematical programs with complementarity constraints, technical report
  version including appendices, arXiv:2009.14047v4, 2021.

\bibitem{kkr23} N. Kreji\'c, E. H. M. Krulikovski, and M. Raydan, A low-cost alternating projection
  approach for a continuous formulation of convex and cardinality constrained optimization,
  {\it Operations Research Forum} 4, article 73, 2023.

\bibitem{kkr25} N. Kreji\'c, E. H. M. Krulikovski, and M. Raydan, An augmented Lagrangian approach
  for cardinality constrained minimization applied to variable selection problems,
  {\it Applied Numerical Mathematics} 208, pp.~284--296, 2025.

\bibitem{krs21} E. H. M. Krulikovski, A. A. Ribeiro, and M. Sachine, On the weak stationarity
  conditions for mathematical programs with cardinality constraints: a unified approach,
  {\it Applied Mathematics \& Optimization} 84, pp.~3451--3473, 2021.

\bibitem{lls21} M. Lapucci, T. Levato, and M. Sciandrone, Convergent inexact penalty decomposition
  methods for cardinality-constrained problems, {\it Journal of Optimization Theory and
  Applications} 188, pp.~473--496, 2021.

\bibitem{llrs23} M. Lapucci, T. Levato, F. Rinaldi, and M. Sciandrone, A unifying framework for
  sparsity-constrained optimization, {\it Journal of Optimization Theory and Applications} 199,
  pp.~663--692, 2023.

\bibitem{lz13} Z. Lu and Y. Zhang, Sparse approximation via penalty decomposition methods,
  {\it SIAM Journal on Optimization} 23, pp.~2448--2478, 2013.

\bibitem{lpr96} Z.-Q. Luo, J.-S. Pang, and D. Ralph, {\it Mathematical Programs with Equilibrium
  Constraints}, Cambridge University Press, Cambridge, 1996.

\bibitem{markowitz} H. Markowitz, Portfolio selection, {\it The Journal of Finance} 7, pp.~77--91,
  1952.

\bibitem{mrss25} J. C. A. Medeiros, A. A. Ribeiro, M. Sachine, and L. D. Secchin, A practical
  second-order optimality condition for cardinality-constrained problems with application to an
  augmented Lagrangian method, {\it Journal of Optimization Theory and Applications} 206,
  article 22, 2025.

\bibitem{mehraprescott} R. Mehra and E. C. Prescott, The equity premium: a puzzle, {\it Journal of
  Monetary Economics} 15, pp.~145--161, 1985.

\bibitem{nl26} A. Nurkanovi\'c and S. Leyffer, A globally convergent method for computing
  B-stationary points of mathematical programs with equilibrium constraints, {\it Computational
  Optimization and Applications}, 2026, DOI 10.1007/s10589-026-00832-x.

\bibitem{rsk22} A. A. Ribeiro, M. Sachine, and E. H. M. Krulikovski, A comparative study of sequential
  optimality conditions for mathematical programs with cardinality constraints,
  {\it Journal of Optimization Theory and Applications} 192, pp.~1067--1083, 2022.

\bibitem{ss00} H. Scheel and S. Scholtes, Mathematical programs with complementarity constraints:
  stationarity, optimality, and sensitivity, {\it Mathematics of Operations Research} 25,
  pp.~1--22, 2000.

\end{thebibliography}
\end{document}